# UNBOUNDED CRITICAL POINTS FOR A CLASS OF LOWER SEMICONTINUOUS FUNCTIONALS

BENEDETTA PELLACCI AND MARCO SQUASSINA

ABSTRACT. In this paper we prove existence and multiplicity results of unbounded critical points for a general class of weakly lower semicontinuous functionals. We will apply a nonsmooth critical point theory developed in [10, 12, 13] and applied in [8, 9, 20] to treat the case of continuous functionals.

## 1. INTRODUCTION

In this paper we prove existence and multiplicity results of unbounded critical points for a class of weakly lower semicontinuous functionals. Let us consider $\Omega$ a bounded open set in $\mathbb{R}^N$ ($N \geqslant 3$) and define the functional $f : H_0^1(\Omega) \to \mathbb{R} \cup \{+\infty\}$ by

$$f(u) = \int_\Omega j(x, u, \nabla u) - \int_\Omega G(x, u),$$

where $j(x, s, \xi) : \Omega \times \mathbb{R} \times \mathbb{R}^N \to \mathbb{R}$ is a measurable function with respect to $x$ for all $(s, \xi) \in \mathbb{R} \times \mathbb{R}^N$, and of class $C^1$ with respect to $(s, \xi)$ for a.e. $x \in \Omega$. We also assume that for almost every $x$ in $\Omega$ and every $s$ in $\mathbb{R}$

(1.1) $\qquad$ the function $\{\xi \mapsto j(x, s, \xi)\}$ is strictly convex.

Moreover, we suppose that there exist a constant $\alpha_0 > 0$ and a positive increasing function $\alpha \in C(\mathbb{R})$ such that the following hypothesis is satisfied for almost every $x \in \Omega$ and for every $(s, \xi) \in \mathbb{R} \times \mathbb{R}^N$

(1.2) $\qquad \alpha_0 |\xi|^2 \leqslant j(x, s, \xi) \leqslant \alpha(|s|)|\xi|^2.$

The functions $j_s(x, s, \xi)$ and $j_\xi(x, s, \xi)$ denote the derivatives of $j(x, s, \xi)$ with respect of the variables $s$ and $\xi$ respectively. Regarding the function $j_s(x, s, \xi)$, we assume that there exist a positive increasing function $\beta \in C(\mathbb{R})$ and a positive constant $R$ such that the following conditions are satisfied almost everywhere in $\Omega$ and for every $\xi \in \mathbb{R}^N$:

(1.3) $\qquad |j_s(x, s, \xi)| \leqslant \beta(|s|)|\xi|^2, \quad$ for every $s$ in $\mathbb{R}$,

(1.4) $\qquad j_s(x, s, \xi) s \geqslant 0, \qquad$ for every $s$ in $\mathbb{R}$ with $|s| \geqslant R.$

*Date*: November 8, 2018.
1991 *Mathematics Subject Classification.* 35J40; 58E05.
Partially supported by Ministero dell'Università e della Ricerca Scientifica e Tecnologica (40% – 1999) and by Gruppo Nazionale per l'Analisi Funzionale e le sue Applicazioni.

Let us notice that from (1.1) and (1.2) it follows that $j_\xi(x, s, \xi)$ satisfies the following growth condition (see Remark 4.1 for more details)

$$(1.5) \qquad |j_\xi(x, s, \xi)| \leqslant 4\alpha(|s|)|\xi|.$$

The function $G(x, s)$ is the primitive with respect to $s$ such that $G(x, 0) = 0$ of a Carathéodory (i.e. measurable with respect to $x$ and continuous with respect to $s$) function $g(x, s)$. We will study two different kinds of problems, according to different nonlinearities $g(x, s)$, that have a main common feature. Indeed, in both cases we cannot expect to find critical points in $L^\infty(\Omega)$. In order to be more precise, let us consider a first model example of nonlinearity and suppose that there exists $p$ such that

$$(1.6) \quad g_1(x, s) = a(x)\mathrm{arctg}\, s + |s|^{p-2}s, \ \ 2 < p < \frac{2N}{N-2}, \ \ a(x) \in L^{\frac{2N}{N+2}}(\Omega), \ a(x) > 0.$$

Notice that from hypotheses (1.2) and (1.6) it follows that $f$ is weakly lower semi-continuous on $H_0^1(\Omega)$. We will also assume that

$$(1.7) \qquad \lim_{|s| \to +\infty} \frac{\alpha(|s|)}{|s|^{p-2}} = 0.$$

Condition (1.7), together with (1.2), allows $f$ to be unbounded from below, so that we cannot look for a global minimum. Moreover, notice that $g(x, s)$ is odd with respect to $s$, so that it would be natural to expect – if $j(x, -s, -\xi) = j(x, s, \xi)$ – the existence of infinitely many solutions as in the semilinear case (see [1]). Unfortunately, we cannot apply any of the classical results of critical point theory, because our functional $f$ is not of class $C^1$ on $H_0^1(\Omega)$. Indeed, notice that $\int_\Omega j(x, v, \nabla v)$ is not differentiable. More precisely, since $j_\xi(x, s, \xi)$ and $j_s(x, s, \xi)$ are not supposed to be bounded with respect to $s$, the terms $j_\xi(x, u, \nabla u) \cdot \nabla v$ and $j_s(x, u, \nabla u)v$ may not be $L^1(\Omega)$ even if $v \in C_0^\infty(\Omega)$. Notice that if $j_s(x, s, \xi)$ and $j_\xi(x, s, \xi)$ were supposed to be bounded with respect to $s$, $f$ would be Gateaux derivable for every $u$ in $H_0^1(\Omega)$ and along any direction $v \in H_0^1(\Omega) \cap L^\infty(\Omega)$ (see [2], [8], [9], [19], [20] for the study of this class of functionals). While in this case, for every $u \in H_0^1(\Omega)$, $f'(u)v$ does not exist even along directions $v \in H_0^1(\Omega) \cap L^\infty(\Omega)$.

In order to deal with the Euler equation of $f$ let us define the following subspace of $H_0^1(\Omega)$ for a fixed $u$ in $H_0^1(\Omega)$.

$$(1.8) \quad W_u = \{v \in H_0^1(\Omega) \,:\, j_\xi(x, u, \nabla u) \cdot \nabla v \in L^1(\Omega), \ \ j_s(x, u, \nabla u)v \in L^1(\Omega)\}.$$

We will see that $W_u$ is dense in $H_0^1(\Omega)$, and we give the following definition of a generalized solution

**Definition 1.1.** *Let $\Lambda \in H^{-1}(\Omega)$ and assume (1.1), (1.2), (1.3). We say that $u$ is a generalized solution of*

$$\begin{cases} -\mathrm{div}\,(j_\xi(x, u, \nabla u)) + j_s(x, u, \nabla u) = \Lambda & \text{in } \Omega, \\ u = 0 & \text{on } \partial\Omega, \end{cases}$$

*if $u \in H_0^1(\Omega)$ and it results*

$$\begin{cases} j_\xi(x, u, \nabla u) \cdot \nabla u \in L^1(\Omega), \ \ j_s(x, u, \nabla u)u \in L^1(\Omega), \\ \displaystyle\int_\Omega j_\xi(x, u, \nabla u) \cdot \nabla v + \int_\Omega j_s(x, u, \nabla u)v = \langle \Lambda, v \rangle \quad \forall\, v \in W_u. \end{cases}$$

We will prove the following result.

**Theorem 1.2.** *Assume conditions (1.1), (1.2), (1.3), (1.4), (1.6), (1.7). Moreover, suppose that there exist $R', \delta > 0$ such that*

(1.9) $$|s| \geq R' \implies pj(x,s,\xi) - j_s(x,s,\xi)s - j_\xi(x,s,\xi) \cdot \xi \geq \delta |\xi|^2$$

*for a.e. $x \in \Omega$ and for all $(s,\xi) \in \mathbb{R} \times \mathbb{R}^N$. Then, if*

$$j(x,-s,-\xi) = j(x,s,\xi),$$

*there exists a sequence $\{u_h\} \subset H^1_0(\Omega)$ of generalized solutions of*

$(P_1)$ $$\begin{cases} -\mathrm{div}\,(j_\xi(x,u,\nabla u)) + j_s(x,u,\nabla u) = g_1(x,u) & \text{in } \Omega, \\ u = 0 & \text{on } \partial\Omega, \end{cases}$$

*with $f(u_h) \to +\infty$.*

In the nonsymmetric case we consider a different class of nonlinearities $g(x,s)$. A simple model example can be the following
(1.10)
$$g_2(x,s) = d(x)\mathrm{arctg}(s^2) + |s|^{p-2}s, \ 2 < p < \frac{2N}{N-2}, \ d(x) \in L^{\frac{N}{2}}(\Omega), \ d(x) > 0.$$

We will prove the following result.

**Theorem 1.3.** *Assume conditions (1.1), (1.2), (1.3), (1.4), (1.7), (1.9), (1.10). Then there exists a nontrivial generalized solution of the following problem*

$(P_2)$ $$\begin{cases} -\mathrm{div}\,(j_\xi(x,u,\nabla u)) + j_s(x,u,\nabla u) = g_2(x,u) & \text{in } \Omega, \\ u = 0 & \text{on } \partial\Omega, \end{cases}$$

Since the functions $\alpha(|s|)$ and $\beta(|s|)$ in (1.2) and (1.3) are not supposed to be bounded, we are dealing with integrands $j(x,s,\xi)$, which may be unbounded with respect to $s$. This class of functionals has also been treated in [3], [4] and [5]. In these papers the existence of a nontrivial solution $u \in L^\infty(\Omega)$ is proved when $g(x,s) = |s|^{p-2}s$. Note that, in this case it is natural to expect solutions in $L^\infty(\Omega)$. In order to prove the existence result, in [4] and [5], a fundamental step is to prove that every cluster point of a Palais-Smale sequence belongs to $L^\infty(\Omega)$. That is, to prove that $u$ is bounded before knowing that it is a solution. Notice that if $u$ is in $L^\infty(\Omega)$ and $v \in C_0^\infty(\Omega)$ then $j_\xi(x,u,\nabla u) \cdot \nabla v$ and $j_s(x,u,\nabla u)v$ are in $L^1(\Omega)$. Therefore, if $g(x,s) = |s|^{p-2}s$, it would be possible to define a solution as a function $u \in L^\infty(\Omega)$ that satisfies the equation associated to $(P_1)$ (or $(P_2)$) in the distributional sense. In our case the function $a(x)$ in (1.6) belongs to $L^{2N/(N+2)}(\Omega)$, so that we can only expect to find solutions in $H^1_0(\Omega)$. In the same way, the function $d(x)$ in (1.10) is in $L^{N/2}(\Omega)$ and also in this case the solutions are not expected to be in $L^\infty(\Omega)$. For these reasons, we have given a definition of solution weaker than the distributional one and we have considered the subspace $W_u$ as the space of the admissible test functions. Notice that if $u \in H^1_0(\Omega)$ is a generalized solution of problem $(P_1)$ (or $(P_2)$) and $u \in L^\infty(\Omega)$, then $u$ is a distributional solution of $(P_1)$ (or $(P_2)$).

We want to stress that we have considered here particular nonlinearities (i.e. $g(x,s) = g_{1,2}(x,s)$) just to present - in a simple case - the main difficulties we are going to tackle. Indeed, Theorems 1.2 and 1.3 will be proved as consequences of two general results (Theorems 2.1 and 2.3). In order to prove these general results we will use an abstract critical point theory for lower semicontinuous functionals developed in [10, 12, 13]. So, first we will show that the functional $f$ can be studied by means of this theory (see Theorem 3.11). Then, we will give a definition of a Palais-Smale sequence $\{u_n\}$ suitable to this situation (Definition 6.3), and we will prove that $u_n$ is compact in $H_0^1(\Omega)$ (Theorems 5.1 and 6.9). In order to do this we will follow the arguments of [8, 9, 19, 20] where the case in which $\alpha(s)$ and $\beta(s)$ are bounded is studied. In our case we will have to modify the test functions used in these papers in order to prove our compactness result. Indeed, here the main difficulty is to find suitable approximations of $u_n$ that belong to $W_{u_n}$, in order to choose them as test functions. For this reason a large amount of work (Theorems 4.7, 4.8, 4.9 and 4.10) is devoted to find possible improvements of the class of allowed test functions. Finally, in Section 7 we will prove a summability result concerning the solutions.

The paper is organized as follows.
In Section 2 we define our general functional, we set the general problem (Problem $(P)$) that we will study and we state the main existence results that we will prove.
In Section 3 we recall (from [10], [12] and [13]) the principal abstract notions and results that we will apply. Moreover, we will study the functional $J : H_0^1(\Omega) \to \mathbb{R} \cup \{+\infty\}$ defined by

$$J(v) = \int_\Omega j(x, v, \nabla v),$$

and we will prove (see Theorem 3.11) that $J$ satisfies a fundamental condition (condition (3.3)) required in order to apply all the abstract results of Section 3.

In Section 4 we find the conditions under which we can compute the directional derivatives of $J$ (Proposition 4.4). Then, we will prove a fundamental inequality regarding the directional derivatives (Proposition 4.5). Moreover, we will prove some Brezis–Browder ([7]) type results (see Theorems 4.7, 4.8, 4.9, 4.10). These results will be important when determining the class of admissible test functions for Problem $(P)$. In particular, in Theorem 4.9 and 4.10 we study the conditions under which we can give a distributional interpretation of Problem $(P)$.

In Section 5 we will prove a compactness result for $J$ (Theorem 5.1). This Theorem will be used to prove that $f$ satisfies our generalized Palais–Smale condition.

In Section 6 we will give the proofs of our general results 2.1 and 2.3 . Then, we will prove Theorems 1.2 and 1.3.

Finally, in Section 7 we will prove a summability result (Theorem 7.1) for a generalized solution in dependence of the summability of the function $g(x,s)$.

## 2. General Setting and Main Results

Let us consider $\Omega$ a bounded open set in $\mathbb{R}^N$ ($N \geqslant 3$). Throughout the paper, we will denote by $\|\cdot\|_p$, $\|\cdot\|_{1,2}$ and $\|\cdot\|_{-1,2}$ the standard norms of the spaces $L^p(\Omega)$, $H_0^1(\Omega)$ and $H^{-1}(\Omega)$ respectively.
Let us define the functional $J : H_0^1(\Omega) \to \mathbb{R} \cup \{+\infty\}$ by

$$(2.1) \qquad J(v) = \int_\Omega j(x, v, \nabla v),$$

where $j(x, s, \xi)$ satisfies hypotheses (1.1), (1.2), (1.3), (1.4). We will prove existence and multiplicity results of generalized solutions (see Definition 1.1) of the problem

$$(P) \qquad \begin{cases} -\operatorname{div}(j_\xi(x, u, \nabla u)) + j_s(x, u, \nabla u) = g(x, u) & \text{in } \Omega, \\ u = 0 & \text{on } \partial\Omega. \end{cases}$$

In order to do this, we will use variational methods, so that we will study the functional $f : H_0^1(\Omega) \to \mathbb{R} \cup \{+\infty\}$ defined by

$$f(v) = J(v) - \int_\Omega G(x,v),$$

where $G(x,s) = \int_0^s g(x,t)dt$ is the primitive of the Carathéodory function $g(x,s)$.

In order to state our multiplicity result let us suppose that $g(x,s)$ satisfies the following conditions.

Assume that for every $\varepsilon > 0$ there exists $a_\varepsilon \in L^{\frac{2N}{N+2}}(\Omega)$ such that

$$(2.2) \qquad |g(x,s)| \leqslant a_\varepsilon(x) + \varepsilon |s|^{\frac{N+2}{N-2}}$$

for a.e. $x \in \Omega$ and every $s \in \mathbb{R}$. Moreover, there exist $p > 2$ and functions $a_0(x), \overline{a}(x) \in L^1(\Omega)$, $b_0(x), \overline{b}(x) \in L^{\frac{2N}{N+2}}(\Omega)$ and $k(x) \in L^\infty(\Omega)$ with $k(x) > 0$ almost everywhere, such that

$$(2.3) \qquad pG(x,s) \leqslant g(x,s)s + a_0(x) + b_0(x)|s|,$$
$$(2.4) \qquad G(x,s) \geqslant k(x)|s|^p - \overline{a}(x) - \overline{b}(x)|s|,$$

for a.e. $x \in \Omega$ and every $s \in \mathbb{R}$ (the constant $p$ is the same of the one in (1.9)).

In this case we will prove the following result.

**Theorem 2.1.** *Assume conditions (1.1), (1.2), (1.3), (1.4), (1.7), (1.9), (2.2), (2.3), (2.4). Moreover, let*

$$(2.5) \qquad j(x,-s,-\xi) = j(x,s,\xi), \qquad g(x,-s) = -g(x,s),$$

*for a.e. $x \in \Omega$ and every $(s,\xi) \in \mathbb{R} \times \mathbb{R}^N$. Then there exists a sequence $\{u_h\} \subset H_0^1(\Omega)$ of generalized solutions of problem (P) with $f(u_h) \to +\infty$.*

**Remark 2.2.** *In the classical results of critical point theory different conditions from (2.2), (2.3) and (2.4) are usually supposed. Indeed, as a growth condition on $g(x,s)$ it is assumed that*

$$(2.6) \quad |g(x,s)| \leqslant a(x) + b|s|^{\sigma-1}, \qquad 2 < \sigma < \frac{2N}{N-2}, \ b \in \mathbb{R}^+, \ a(x) \in L^{\frac{2N}{N+2}}(\Omega).$$

*Note that (2.6) implies (2.2). Indeed, suppose that $g(x,s)$ satisfies (2.6), then Young inequality implies that (2.2) is satisfied with $a_\varepsilon(x) = a(x) + C(b,\varepsilon)$.*

*Moreover, as a superlinear condition it is usually assumed that there exist $p > 2$ and $R > 0$ such that*

$$(2.7) \qquad 0 < pG(x,s) \leqslant g(x,s)s, \qquad \text{for every } s, \text{ with } |s \geqslant R.$$

*Note that this condition is stronger than conditions (2.3), (2.4). Indeed, suppose that $g(x,s)$ satisfies (2.7) and notice that this implies that there exists $a_0 \in L^1(\Omega)$ such that*

$$pG(x,s) \leqslant g(x,s)s + a_0(x), \qquad \text{for every } s \text{ in } \mathbb{R}.$$

*Then (2.3) is satisfied with $b_0(x) \equiv 0$. Moreover, from (2.7) we deduce that there exists $\overline{a}(x) \in L^1(\Omega)$ such that*

$$G(x,s) \geqslant \frac{1}{R^p} \min\{G(x,R), G(x,-R), 1\}|s|^p - \overline{a}(x).$$

*then also (2.4) is satisfied.*

In order to state our existence result in the nonsymmetric case, assume that the function $g$ satisfies the following condition

$$|g(x,s)| \leqslant a_1(x)|s| + b|s|^{\sigma-1}, \tag{2.8}$$

$$2 < \sigma < \frac{2N}{N-2}, \quad a_1(x) \in L^{\frac{N}{2}}(\Omega), \ b \in \mathbb{R}^+.$$

We will prove the following result.

**Theorem 2.3.** *Assume conditions* (1.1), (1.2), (1.3), (1.4), (1.7), (1.9), (2.3), (2.4), (2.8). *Moreover, let*

$$\lim_{s \to 0} \frac{g(x,s)}{s} = 0 \quad a.e. \ in \ \Omega. \tag{2.9}$$

*Then there exists a nontrivial generalized solution of the problem* (P).
*In addition, there exist* $\varepsilon > 0$ *such that for every* $\Lambda \in H^{-1}(\Omega)$ *with* $\|\Lambda\|_{-1,2} < \varepsilon$ *the problem*

$$(P_\Lambda) \quad \begin{cases} -\operatorname{div}(j_\xi(x,u,\nabla u)) + j_s(x,u,\nabla u) = g(x,u) + \Lambda & \text{in } \Omega, \\ u = 0 & \text{on } \partial\Omega, \end{cases}$$

*has at least two generalized solutions* $u_1, u_2$ *with* $f(u_1) \leqslant 0 < f(u_2)$.

**Remark 2.4.** *Notice that, in order to have* $g(x,v)v \in L^1(\Omega)$ *for every* $v \in H^1_0(\Omega)$, *the function* $a_1(x)$ *has to be in* $L^{\frac{N}{2}}(\Omega)$. *Nevertheless, also in this case we cannot expect to find bounded solution of problem* (P). *The situation is even worse in problem* $(P_\Lambda)$, *indeed in this case we can only expect to find solutions that belong to* $H^1_0(\Omega) \cap \operatorname{dom}(J)$.

**Remark 2.5.** *Notice that condition* (2.8) *implies* (2.2). *Indeed, suppose that* $g(x,s)$ *satisfies* (2.8). *Then Young inequality implies that for every* $\varepsilon > 0$ *we have*

$$|g(x,s)| \leqslant \beta(\varepsilon)(a_1(x))^{\frac{N+2}{4}} + \varepsilon|s|^{\frac{N+2}{N-2}} + \gamma(\varepsilon,b),$$

*where* $\beta(\varepsilon)$ *and* $\gamma(\varepsilon,b)$ *are positive constants depending on* $\varepsilon, b$. *Now, since* $a_1(x) \in L^{\frac{N}{2}}(\Omega)$, *it results*

$$a_\varepsilon(x) := \left(\beta(\varepsilon)(a_1(x))^{\frac{N+2}{4}} + \gamma(\varepsilon,b)\right) \in L^{\frac{2N}{N+2}}(\Omega),$$

*so* (2.2) *holds.*

## 3. Abstract results of critical point theory

In this section we will recall the principal abstract notions and results that we will use in the sequel. We refer to [10], [12] and [13], where this abstract theory is developed.

Moreover, we will prove that our functional $f$ satisfies a fundamental condition (see condition (3.3) and Theorem 3.11) requested to apply the abstract results.

Let $X$ be a metric space and let $f : X \to \mathbb{R} \cup \{+\infty\}$ be a lower semicontinuous function. We set

$$\operatorname{dom}(f) = \{u \in X : f(u) < +\infty\}, \quad \operatorname{epi}(f) = \{(u,\eta) \in X \times \mathbb{R} : f(u) \leqslant \eta\}.$$

The set $\operatorname{epi}(f)$ is endowed with the metric

$$d((u,\eta),(v,\mu)) = \left(d(u,v)^2 + (\eta-\mu)^2\right)^{1/2}.$$

Let us define the function $\mathcal{G}_f : \operatorname{epi}(f) \to \mathbb{R}$ by

$$\mathcal{G}_f(u, \eta) := \eta. \tag{3.1}$$

Note that $\mathcal{G}_f$ is Lipschitz continuous of constant 1.

From now on we denote with $B(u, \delta)$ the open ball of center $u$ and of radius $\delta$.

Let us recall the definition of the weak slope for a continuous function introduced in [10, 12, 15, 16].

**Definition 3.1.** *Let $X$ be a complete metric space, $g : X \to \mathbb{R}$ be a continuous function, and $u \in X$. We denote by $|dg|(u)$ the supremum of the real numbers $\sigma$ in $[0, +\infty)$ such that there exist $\delta > 0$ and a continuous map*

$$\mathcal{H} : B(u, \delta) \times [0, \delta] \to X,$$

*such that, for every $v$ in $B(u, \delta)$, and for every $t$ in $[0, \delta]$ it results*

$$d(\mathcal{H}(v, t), v) \leqslant t,$$
$$g(\mathcal{H}(v, t)) \leqslant g(v) - \sigma t.$$

*The extended real number $|dg|(u)$ is called the weak slope of $g$ at $u$.*

According to the previous definition, for every lower semicontinuous function $f$ we can consider the metric space $\operatorname{epi}(f)$ so that the weak slope of $\mathcal{G}_f$ is well defined. Therefore, we can define the weak slope of a lower semicontinuous function $f$ by using $|d\mathcal{G}_f|(u, f(u))$. More precisely, we have the following definition.

**Definition 3.2.** *For every $u \in \operatorname{dom}(f)$ let*

$$|df|(u) = \begin{cases} \dfrac{|d\mathcal{G}_f|(u, f(u))}{\sqrt{1 - |d\mathcal{G}_f|(u, f(u))^2}} & \text{if } |d\mathcal{G}_f|(u, f(u)) < 1, \\ +\infty & \text{if } |d\mathcal{G}_f|(u, f(u)) = 1. \end{cases}$$

The previous notions allow us to give the following definitions.

**Definition 3.3.** *Let $X$ be a complete metric space and $f : X \to \mathbb{R} \cup \{+\infty\}$ a lower semicontinuous function. We say that $u \in \operatorname{dom}(f)$ is a (lower) critical point of $f$ if $|df|(u) = 0$. We say that $c \in \mathbb{R}$ is a (lower) critical value of $f$ if there exists a (lower) critical point $u \in \operatorname{dom}(f)$ of $f$ with $f(u) = c$.*

**Definition 3.4.** *Let $X$ be a complete metric space, $f : X \to \mathbb{R} \cup \{+\infty\}$ a lower semicontinuous function and let $c \in \mathbb{R}$. We say that $f$ satisfies the Palais–Smale condition at level $c$ ($(PS)_c$ in short), if every sequence $\{u_n\}$ in $\operatorname{dom}(f)$ such that*

$$|df|(u_n) \to 0,$$
$$f(u_n) \to c,$$

*admits a subsequence $\{u_{n_k}\}$ converging in $X$.*

For every $\eta \in \mathbb{R}$, let us define the set

$$f^\eta = \{u \in X : f(u) < \eta\}. \tag{3.2}$$

The following result gives a criterion to obtain a lower estimate of $|df|(u)$. For the proof see [12].

**Proposition 3.5.** *Let $f : X \to \mathbb{R} \cup \{+\infty\}$ be a lower semicontinuous function defined on a complete metric space, and let $u \in \mathrm{dom}(f)$. Let us assume that there exist $\delta > 0$, $\eta > f(u)$, $\sigma > 0$ and a continuous function $\mathcal{H} : B(u,\delta) \cap f^\eta \times [0,\delta] \to X$ such that*

$$d(\mathcal{H}(v,t),v) \leqslant t, \qquad \forall\, v \in B(u,\delta) \cap f^\eta,$$
$$f(\mathcal{H}(v,t)) \leqslant f(v) - \sigma t, \qquad \forall\, v \in B(u,\delta) \cap f^\eta.$$

*Then $|df|(u) \geqslant \sigma$.*

We will also use the notion of equivariant weak slope (see [9]).

**Definition 3.6.** *Let $X$ be a normed linear space and $f : X \to \mathbb{R} \cup \{+\infty\}$ an even lower semicontinuous function with $f(0) < +\infty$. For every $(0,\eta) \in \mathrm{epi}\,(f)$ we denote with $|d_{\mathbb{Z}_2}\mathcal{G}_f|(0,\eta)$ the supremum of the numbers $\sigma \in [0,\infty)$ such that there exist $\delta > 0$ and a continuous map*

$$\mathcal{H} = (\mathcal{H}_1, \mathcal{H}_2) : (B((0,\eta),\delta) \cap \mathrm{epi}\,(f)) \times [0,\delta] \to \mathrm{epi}\,(f)$$

*satisfying*

$$d(\mathcal{H}((w,\mu),t),(w,\mu)) \leqslant t, \qquad \mathcal{H}_2((w,\mu),t) \leqslant \mu - \sigma t,$$
$$\mathcal{H}_1((-w,\mu),t) = -\mathcal{H}_1((w,\mu),t),$$

*for every $(w,\mu) \in B((0,\eta),\delta) \cap \mathrm{epi}\,(f)$ and $t \in [0,\delta]$.*

In order to compute $|d\mathcal{G}_f|(u,\eta)$ it will be useful the following result (for the proof see [12]).

**Proposition 3.7.** *Let $X$ be a normed linear space, $J : X \to \mathbb{R} \cup \{+\infty\}$ a lower semicontinuous functional, $I : X \to \mathbb{R}$ a $C^1$ functional and let $f = J + I$. Then the following facts hold:*

*(a) for every $(u,\eta) \in \mathrm{epi}(f)$ we have*

$$|d\mathcal{G}_f|(u,\eta) = 1 \iff |d\mathcal{G}_J|(u,\eta - I(u)) = 1\,;$$

*(b) if $J$ and $I$ are even, for every $\eta \geqslant f(0)$ we have*

$$|d_{\mathbb{Z}_2}\mathcal{G}_f|(0,\eta) = 1 \iff |d_{\mathbb{Z}_2}\mathcal{G}_J|(0,\eta - I(0)) = 1\,;$$

*(c) if $u \in \mathrm{dom}(f)$ and $I'(u) = 0$, then*

$$|df|(u) = |dJ|(u).$$

*Proof.* Assertions $(a)$ and $(c)$ follow by arguing as in [12]. Assertion $(b)$ can be reduced to $(a)$ after observing that, since $I$ is even, it results $I'(0) = 0$. □

In [10, 12] variational methods for lower semicontinuous functionals are developed. Moreover, it is shown that the following condition is fundamental in order to apply this abstract theory to the study of lower semicontinuous functions.

(3.3) $\qquad \forall (u,\eta) \in \mathrm{epi}\,(f):\ f(u) < \eta \implies |d\mathcal{G}_f|(u,\eta) = 1.$

In the next section we will prove that our functional $f$ satisfies (3.3).

The next result gives a criterion to verify condition (3.3).

**Theorem 3.8.** *Let $(u,\eta) \in \mathrm{epi}(f)$ with $f(u) < \eta$. Assume that for every $\varrho > 0$ there exist $\delta > 0$ and a continuous map*

$$\mathcal{H} : \bigl\{w \in B(u,\delta) :\ f(w) < \eta + \delta\bigr\} \times [0,\delta] \to X$$

*satisfying*

$$d(\mathcal{H}(w,t), w) \leqslant \varrho t, \quad f(\mathcal{H}(w,t)) \leqslant (1-t)f(w) + t(f(u) + \varrho)$$

*whenever* $w \in B(u, \delta)$, $f(w) < \eta + \delta$ *and* $t \in [0, \delta]$.

Then we have $|d\mathcal{G}_f|(u, \eta) = 1$. If moreover $X$ is a normed space, $f$ is even, $u = 0$ and $\mathcal{H}(-w, t) = -\mathcal{H}(w, t)$, then we have $|d_{\mathbb{Z}_2}\mathcal{G}_f|(0, \eta) = 1$

*Proof.* See [13, Corollary 2.11]. □

Let us now recall from [10] the following existence result.

**Theorem 3.9.** *Let $X$ be a complete metric space, $f : X \to \mathbb{R} \cup \{+\infty\}$ a lower semicontinuous function satisfying (3.3).*

*First, suppose that there exist $w \in X$, $\eta > 0$ and $r > 0$ such that*

(3.4) $$f(u) > \eta, \quad \forall u \in X, \ \|u\| = r,$$

(3.5) $$f(w) < \eta, \quad \|w\| > r.$$

*In this case we set*

$$\Gamma = \{\gamma : [0, 1] \to X, \ continuous \ \gamma(0) = 0, \ \gamma(1) = w\}.$$

*Finally, suppose that $f$ satisfies the Palais-Smale condition at the level*

$$c = \inf_{\Gamma} \sup_{[0,1]} f(\gamma(t)) < +\infty,$$

*then, there exists a nontrivial critical point $u$ of $f$ such that $f(u) = c$.*

*Now, assume that there exist $v_0, v_1$ in $X$ and $r > 0$ such that*

$$\|v_0\| < r, \quad \|v_1\| > r, \quad \inf f(\overline{B_r(0)}) > -\infty, \quad and$$

$$\inf\{f(u) : u \in X, \ \|u\| = r\} > \max\{f(v_0), f(v_1)\}.$$

*Let*

$$\Gamma = \left\{\gamma : [0, 1] \to D(f) \ continuous \ with \ \gamma(0) = v_0, \gamma(1) = v_1\right\}$$

*and assume that $\Gamma \neq \emptyset$ and that $f$ satisfies the Palais–Smale condition at the two levels*

$$c_1 = \inf f(\overline{B_r(0)}), \quad c_2 = \inf_{\gamma \in \Gamma} \sup_{[0,1]}(f \circ \gamma) < +\infty.$$

*Then $c_1 < c_2$ and there exist at least two critical points $u_1, u_2$ of $f$ such that $f(u_i) = c_i$ for every $i = 1, 2$.*

*Proof.* The existence of $u_1$ follows by Ekeland Variational Principle, while $u_2$ can be found by applying [10, Theorem 4.5]. □

In the equivariant case we apply the following existence result (see [10, 17]).

**Theorem 3.10.** *Let $X$ be a Banach space and let $f : X \to \mathbb{R} \cup \{+\infty\}$ an even lower semicontinuous function. Let us assume that*

(a) *there exist $\rho > 0$, $\gamma > f(0)$ and a subspace $V \subset X$ of finite codimension such that*

$$\forall u \in V : \|u\| = \rho \implies f(u) \geqslant \gamma;$$

(b) *for every finite dimensional subspace $W \subset X$, there exists $R > 0$ such that*

$$\forall u \in W : \|u\| \geqslant R \implies f(u) \leqslant f(0);$$

(c) $f$ satisfies $(PS)_c$ for any $c \geqslant \gamma$ and $f$ satisfies (3.3);
(d) it results $|d_{\mathbb{Z}_2}\mathcal{G}_f|(0,\eta) \neq 0$ for every $\eta > f(0)$.
Then there exists a sequence $\{u_h\}$ of critical points of $f$ with
$$\lim_{h\to\infty} f(u_h) = +\infty.$$

Now, let us set $X = H_0^1(\Omega)$ and consider the functional $J : H_0^1(\Omega) \to \mathbb{R} \cup \{+\infty\}$ defined in (2.1). From hypothesis (1.2) we obtain that $J$ is lower semicontinuous. We will prove that $J$ satisfies (3.3). In order to do this, it will be useful the following function. For $k \geqslant 1$, let $T_k : \mathbb{R} \to \mathbb{R}$ be the truncation at height $k$, defined by

$$(3.6) \qquad T_k(s) = s \quad \text{if } |s| \leqslant k, \qquad T_k(s) = k\frac{s}{|s|} \quad \text{if } |s| \geqslant k.$$

We will prove the following result.

**Theorem 3.11.** *Assume conditions (1.1), (1.2), (1.4). Then, for every $(u,\eta) \in \mathrm{epi}\,(J)$ with $J(u) < \eta$ it results*
$$|d\mathcal{G}_J|(u,\eta) = 1.$$
*Moreover, if $j(x,-s,-\xi) = j(x,s,\xi)$, $\forall \eta > J(0)(=0)$ it results $|d_{\mathbb{Z}_2}\mathcal{G}_J|(0,\eta) = 1$.*

*Proof.* Let $(u,\eta) \in \mathrm{epi}\,(J)$ with $J(u) < \eta$ and let $\varrho > 0$. Then, there exists $\delta \in (0,1]$, $\delta = \delta(\varrho)$ and $k \geqslant 1$, $k = k(\varrho)$ such that $k \geqslant R$ ($R$ is defined in (1.4)) and
$$(3.7) \qquad \|T_k(v) - v\|_{1,2} < \varrho, \quad \text{for every } v \in B(u,\delta)$$
From (1.2) we have
$$j(x,v,\nabla T_k(v)) \leqslant \alpha(k)|\nabla v|^2.$$
Then, up to reducing $\delta$, we get that the following inequalities hold
$$(3.8) \qquad \int_\Omega j(x,v,\nabla T_k(v)) < \int_\Omega j(x,u,\nabla T_k(u)) + \varrho \leqslant \int_\Omega j(x,u,\nabla u) + \varrho,$$
for every $v \in B(u,\delta)$. Now, let us prove that the following inequality is satisfied for every $t \in [0,\delta]$ and for every $v \in B(u,\delta)$
$$(3.9) \qquad J((1-t)v + tT_k(v)) \leqslant (1-t)J(v) + t(J(u) + \varrho).$$
Since $j(x,s,\xi)$ is of class $C^1$ with respect to the variable $s$ and from (1.1) we obtain that there exists $\theta \in [0,1]$ such that
$$\begin{aligned}
&j(x,(1-t)v + tT_k(v), (1-t)\nabla v + t\nabla T_k(v)) - j(x,v,\nabla v) = \\
&= j(x,(1-t)v + tT_k(v), (1-t)\nabla v + t\nabla T_k(v)) - j(x,v,(1-t)\nabla v + t\nabla T_k(v)) \\
&\quad + j(x,v,(1-t)\nabla v + t\nabla T_k(v)) - j(x,v,\nabla v) \\
&\leqslant tj_s(x,v + \theta t(T_k(v) - v), (1-t)\nabla v + t\nabla T_k(v))(T_k(v) - v) \\
&\quad + t\left(j(x,v,\nabla T_k(v)) - j(x,v,\nabla v)\right).
\end{aligned}$$
Notice that
$$\begin{aligned}
v(x) \geqslant k &\implies v(x) + \theta t(T_k(v(x)) - v(x)) \geqslant k \geqslant R, \\
v(x) \leqslant -k &\implies v(x) + \theta t(T_k(v(x)) - v(x)) \leqslant -k \leqslant -R.
\end{aligned}$$
Then, in view of (1.4) one has
$$j_s(x,v + \theta t(T_k(v) - v), (1-t)\nabla v + t\nabla T_k(v))(T_k(v) - v) \leqslant 0,$$

whence
$$j(x, (1-t)v + tT_k(v), (1-t)\nabla v + t\nabla T_k(v)) \leqslant (1-t)j(x, v, \nabla v) + tj(x, v, \nabla T_k(v)).$$

Therefore from (3.8) one gets (3.9). In order to apply Theorem 3.8 we define
$$\mathcal{H} : \{v \in B(u,\delta) : J(v) < \eta + \delta\} \times [0,\delta] \to H_0^1(\Omega)$$
by setting
$$\mathcal{H}(v,t) = (1-t)v + tT_k(v).$$
Then, taking into account (3.7) and (3.9), it results
$$d(\mathcal{H}(v,t), v) \leqslant \varrho t, \quad J(\mathcal{H}(v,t)) \leqslant (1-t)J(v) + t(J(u) + \varrho)$$
whenever $v \in B(u,\delta)$, $J(v) < \eta + \delta$ and $t \in [0,\delta]$. The first assertion now follows from Theorem 3.8. Finally, since $\mathcal{H}(-v,t) = \mathcal{H}(v,t)$ one also has $|d_{\mathbb{Z}_2}\mathcal{G}_J|(0,\eta) = 1$, whenever $j(x,-s,-\xi) = j(x,s,\xi)$. $\square$

## 4. THE VARIATIONAL SETTING

This section regards the relations between $|dJ|(u)$ and the derivative of the functional $J$. Moreover, we will prove some Brezis-Browder (see [7]) type results.
First of all let us make some remarks.

**Remark 4.1.** *It is readily seen that hypothesis (1.1) and the right inequality of (1.2) imply that there exists a positive increasing function $\overline{\alpha}(|s|)$ such that*

(4.1) $$|j_\xi(x,s,\xi)| \leqslant \overline{\alpha}(|s|)|\xi|$$

*for a.e. $x \in \Omega$ and every $(s,\xi) \in \mathbb{R} \times \mathbb{R}^N$. Indeed, from (1.1) one has*
$$\forall v \in \mathbb{R}^N : |v| \leqslant 1 \implies j(x,s,\xi + |\xi|v) \geqslant j(x,s,\xi) + j_\xi(x,s,\xi) \cdot v|\xi|.$$
*This, and (1.2) yield*
$$j_\xi(x,s,\xi) \cdot v|\xi| \leqslant 4\alpha(|s|)|\xi|^2.$$
*From the arbitrariness of $v$ (4.1) follows. On the other hand, if (4.1) holds we have*
$$|j(x,s,\xi)| \leqslant \int_0^1 |j_\xi(x,s,t\xi) \cdot \xi| dt \leqslant \frac{1}{2}\overline{\alpha}(|s|)|\xi|^2.$$
*So that $j(x,s,\xi)$ satisfies a growth condition equivalent to (4.1). For these reasons it is not restrictive to suppose that the functions in the right hand side of (1.2) and (4.1) are the same. Note that, in particular, $j_\xi(x,s,0) = 0$.*

**Remark 4.2.** *Notice that it is not restrictive to suppose that the functions $\alpha(s)$ and $\beta(s)$ are increasing. Indeed, if this is not the case we can consider the functions*
$$A_r(|s|) = \sup_{|s|\leqslant r} \alpha(|s|), \quad B_r(|s|) = \sup_{|s|\leqslant r} \beta(|s|),$$
*which are increasing.*

**Remark 4.3.** *Let us point out that the assumption of strict convexity on the function $\{\xi \to j(x,s,\xi)\}$ implies that for almost every $x$ in $\Omega$ and for every $s$ in $\mathbb{R}$, we have:*

(4.2) $$[j_\xi(x,s,\xi) - j_\xi(x,s,\xi^*)] \cdot (\xi - \xi^*) > 0,$$

*for every $\xi, \xi^* \in \mathbb{R}^N$, with $\xi \neq \xi^*$. Moreover, hypotheses (1.1) and (1.2) imply that,*

(4.3) $$j_\xi(x,s,\xi) \cdot \xi \geqslant \alpha_0|\xi|^2.$$

*Indeed, we have*

$$0 = j(x, s, 0) \geq j(x, s, \xi) + j_\xi(x, s, \xi) \cdot (0 - \xi),$$

*then by hypothesis (1.2) we get (4.3).*

Now, let us set for every $u \in H^1_0(\Omega)$ the subspace

(4.4) $\qquad V_u = \left\{ v \in H^1_0(\Omega) \cap L^\infty(\Omega) : u \in L^\infty(\{x \in \Omega : v(x) \neq 0\}) \right\}.$

As proved in [14], $V_u$ is a vector space dense in $H^1_0(\Omega)$. Since $V_u \subset W_u$, also $W_u$ (see the Introduction) is dense in $H^1_0(\Omega)$. In the following proposition we study the conditions under which we can compute the directional derivative of $J$.

**Proposition 4.4.** *Assume conditions (1.2), (1.3), (1.5). Then there exists $J'(u)(v)$ for every $u \in \mathrm{dom}(J)$, $v \in V_u$ and one has $j_s(x, u, \nabla u)v \in L^1(\Omega)$, $j_\xi(x, u, \nabla u) \cdot \nabla v \in L^1(\Omega)$ and*

$$J'(u)(v) = \int_\Omega j_\xi(x, u, \nabla u) \cdot \nabla v + \int_\Omega j_s(x, u, \nabla u)v.$$

*Proof.* Let $u \in \mathrm{dom}(J)$ and $v \in V_u$. For every $t \in \mathbb{R}$ and a.e. $x \in \Omega$, we set

$$F(x, t) = j(x, u(x) + tv(x), \nabla u(x) + t\nabla v(x)).$$

Since $v \in V_u$ and by using (1.2), it follows that $F(x, t) \in L^1(\Omega)$. Moreover, it results

$$\frac{\partial F}{\partial t}(x, t) = j_s(x, u + tv, \nabla u + t\nabla v)v + j_\xi(x, u + tv, \nabla u + t\nabla v) \cdot \nabla v.$$

From hypotheses (1.3) and (1.5) we get that for every $x \in \Omega$ with $v(x) \neq 0$, it results

$$\left| \frac{\partial F}{\partial t}(x, t) \right| \leq \|v\|_\infty \beta(\|u\|_\infty + \|v\|_\infty)(|\nabla u| + |\nabla v|)^2$$
$$+ \alpha(\|u\|_\infty + \|v\|_\infty)(|\nabla u| + |\nabla v|)|\nabla v|.$$

Since the function in the right hand side of the previous inequality is in $L^1(\Omega)$, the assertion follows. $\square$

In the sequel we will often use the function $H \in C^\infty(\mathbb{R})$ given by

(4.5) $\qquad H(s) = 1, \quad \text{on } [-1, 1], \quad H(s) = 0, \text{ outside } [-2, 2] \text{ and } |H'(s)| \leq 2.$

Now we can prove a fundamental inequality regarding $|dJ|(u)$.

**Proposition 4.5.** *Assume conditions (1.2), (1.3), (1.5). Then we have*

$$|d(J - w)|(u) \geq$$
$$\geq \sup \left\{ \int_\Omega j_\xi(x, u, \nabla u) \cdot \nabla v + \int_\Omega j_s(x, u, \nabla u)v - \langle w, v \rangle : v \in V_u, \|v\|_{1,2} \leq 1 \right\}$$

*for every $u \in \mathrm{dom}(J)$ and every $w \in H^{-1}(\Omega)$.*

*Proof.* If it results $|d(J - w)|(u) = +\infty$, or if it holds

$$\sup \left\{ \int_\Omega j_\xi(x, u, \nabla u) \cdot \nabla v + \int_\Omega j_s(x, u, \nabla u)v - \langle w, v \rangle : v \in V_u, \|v\|_{1,2} \leq 1 \right\} = 0,$$

the inequality holds. Otherwise, let $u \in \mathrm{dom}(J)$ and let $\eta \in \mathbb{R}^+$ be such that $J(u) < \eta$. Since we want to apply Proposition 3.5, let us consider $J^\eta$ defined in (3.2). Moreover, let us consider $\overline{\sigma} > 0$ and $\overline{v} \in V_u$ such that $\|\overline{v}\|_{1,2} \leq 1$ and

(4.6) $\qquad \int_\Omega j_\xi(x, u, \nabla u) \cdot \nabla \overline{v} + \int_\Omega j_s(x, u, \nabla u)\overline{v} - \langle w, \overline{v} \rangle < -\overline{\sigma}.$

Let us fix $\varepsilon > 0$ and let us prove that there exists $k_0 \in \mathbb{N}$ such that

(4.7) $$\left\| H\left(\frac{u}{k_0}\right) \overline{v} \right\|_{1,2} < 1 + \varepsilon,$$

(4.8) $$\int_\Omega j_s(x, u, \nabla u) H\left(\frac{u}{k_0}\right) \overline{v}$$
$$+ \int_\Omega j_\xi(x, u, \nabla u) \cdot \nabla \left( H\left(\frac{u}{k_0}\right) \overline{v} \right) - \left\langle w, H\left(\frac{u}{k_0}\right) \overline{v} \right\rangle < -\overline{\sigma}.$$

Let us set $v_k = H\left(\frac{u}{k}\right) \overline{v}$, where $H(s)$ is defined in (4.5). Since $\overline{v} \in V_u$ we deduce that $v_k \in V_u$ for every $k \geqslant 1$ and $v_k$ converges to $\overline{v}$ in $H_0^1(\Omega)$. This, together with the fact that $\|\overline{v}\|_{1,2} \leqslant 1$ implies (4.7). Moreover, Proposition 4.4 implies that we can consider $J'(u)(v_k)$. In addition, we have, as $k$ tends to infinity

$$j_s(x, u(x), \nabla u(x)) v_k(x) \to j_s(x, u(x), \nabla u(x)) \overline{v}(x), \quad \text{a.e. in } \Omega,$$
$$j_\xi(x, u(x), \nabla u(x)) \cdot \nabla v_k(x) \to j_\xi(x, u(x), \nabla u(x)) \cdot \nabla \overline{v}(x), \quad \text{a.e. in } \Omega.$$

Moreover, we get

$$\left| j_s(x, u, \nabla u) H\left(\frac{u}{k}\right) \overline{v} \right| \leqslant |j_s(x, u, \nabla u) \overline{v}|$$
$$|j_\xi(x, u, \nabla u) \nabla v_k| \leqslant |j_\xi(x, u, \nabla u)| |\nabla \overline{v}| + 2|\overline{v}| |j_\xi(x, u, \nabla u) \cdot \nabla u|$$

Since $\overline{v} \in V_u$ and by using (1.3) and (1.5), we can apply Lebesgue Dominated Convergence Theorem to obtain

$$\lim_{k \to +\infty} \int_\Omega j_s(x, u, \nabla u) u_k = \int_\Omega j_s(x, u, \nabla u) \overline{v},$$
$$\lim_{k \to +\infty} \int_\Omega j_\xi(x, u, \nabla u) \cdot \nabla u_k = \int_\Omega j_\xi(x, u, \nabla u) \cdot \nabla \overline{v},$$

which, together with (4.6), implies (4.8). Let us now show that there exists $\delta_1 > 0$ such that

(4.9) $$\left\| H\left(\frac{z}{k_0}\right) \overline{v} \right\| \leqslant 1 + \varepsilon,$$

(4.10) $$\int_\Omega j_\xi(x, z, \nabla z) \cdot \nabla \left( H\left(\frac{z}{k_0}\right) \overline{v} \right)$$
$$+ \int_\Omega j_s(x, z, \nabla z) H\left(\frac{z}{k_0}\right) \overline{v} - \left\langle w, H\left(\frac{z}{k_0}\right) \overline{v} \right\rangle < -\overline{\sigma},$$

for every $z \in B(u, \delta_1) \cap J^\eta$. Indeed, take $u_n \in J^\eta$ such that $u_n \to u$ in $H_0^1(\Omega)$ and set

$$v_n = H\left(\frac{u_n}{k_0}\right) \overline{v}.$$

We have that $v_n \to H\left(\frac{u}{k_0}\right) \overline{v}$ in $H_0^1(\Omega)$, so that (4.9) follows from (4.7). Moreover, note that $v_n \in V_{u_n}$, so that from Proposition 4.4 we deduce that we can consider $J'(u_n)(v_n)$. From (1.3) and (1.5) it follows

$$|j_s(x, u_n, \nabla u_n) v_n| \leqslant \beta(2k_0) \|\overline{v}\|_\infty |\nabla u_n|^2,$$
$$|j_\xi(x, u_n, \nabla u_n) \cdot \nabla v_n| \leqslant \alpha(2k_0) |\nabla u_n| \left[ \frac{2}{k_0} \|\overline{v}\|_\infty |\nabla u_n| + |\nabla \overline{v}| \right].$$

Then
$$\lim_{n\to+\infty}\int_\Omega j_s(x,u_n,\nabla u_n)v_n = \int_\Omega j_s(x,u,\nabla u)H\left(\frac{u}{k_0}\right)\overline{v},$$
$$\lim_{n\to+\infty}\int_\Omega j_\xi(x,u_n,\nabla u_n)\cdot\nabla v_n = \int_\Omega j_\xi(x,u,\nabla u)\cdot\nabla\left[H\left(\frac{u}{k_0}\right)\overline{v}\right],$$

which, together with (4.8), implies (4.10). Now, observe that (4.10) is equivalent to say that $J'(z)\left(H\left(\frac{z}{k}\right)\overline{v}\right) - \langle w, H\left(\frac{z}{k}\right)\overline{v}\rangle < -\overline{\sigma}$. Thus, there exists $\delta < \delta_1$ with

$$(4.11)\qquad J\left(z + \frac{t}{1+\varepsilon}H\left(\frac{z}{k_0}\right)\overline{v}\right) - J(z) - \langle w, \frac{t}{1+\varepsilon}H\left(\frac{z}{k}\right)\overline{v}\rangle \leqslant -\frac{\overline{\sigma}}{1+\varepsilon}t,$$

for every $t \in [0,\delta]$ and $z \in B(u,\delta) \cap J^\eta$. Finally, let us define the continuous function $\mathcal{H}: B(u,\delta) \cap J^\eta \times [0,\delta] \to H_0^1(\Omega)$ by

$$\mathcal{H}(z,t) = z + \frac{t}{1+\varepsilon}H\left(\frac{z}{k_0}\right)\overline{v}.$$

From (4.9) and (4.11) we deduce that $\mathcal{H}$ satisfies all the hypotheses of Proposition 3.5. Then, $|dJ|(u) > \frac{\sigma}{1+\varepsilon}$, and the conclusion follows from the arbitrariness of $\varepsilon$. □

The next Lemma will be useful in proving two Brezis–Browder (see [7]) type results for $J$.

**Lemma 4.6.** *Assume conditions (1.1), (1.2), (1.3), (1.4) and let $u \in \mathrm{dom}(J)$. Then*

$$(4.12)\qquad \int_\Omega j_\xi(x,u,\nabla u)\cdot\nabla u + \int_\Omega j_s(x,u,\nabla u)u \leqslant |dJ|(u)\|u\|_{1,2}.$$

*In particular, if $|dJ|(u) < +\infty$ it results*
$$j_\xi(x,u,\nabla u)\cdot\nabla u \in L^1(\Omega),\quad j_s(x,u,\nabla u)u \in L^1(\Omega).$$

*Proof.* First, notice that if $u$ is such that $|dJ|(u) = +\infty$, or
$$\int_\Omega j_\xi(x,u,\nabla u)\cdot\nabla u + \int_\Omega j_s(x,u,\nabla u)u \leqslant 0,$$
then the conclusion holds. Otherwise, let us consider $k \geqslant 1$, $u \in \mathrm{dom}(J)$ with $|dJ|(u) < +\infty$, and $\sigma > 0$ such that
$$\int_\Omega j_\xi(x,u,\nabla u)\cdot\nabla T_k(u) + \int_\Omega j_s(x,u,\nabla u)T_k(u) > \sigma\|T_k(u)\|_{1,2},$$
where $T_k(s)$ is defined in (3.6). We will prove that $|dJ|(u) \geqslant \sigma$.
Let us fix $\varepsilon > 0$, we first want to show that there exists $\delta_1 > 0$ such that

$$(4.13)\qquad \|T_k(w)\|_{1,2} \leqslant (1+\varepsilon)\|T_k(u)\|_{1,2}$$

$$(4.14)\qquad \int_\Omega j_\xi(x,w,\nabla w)\cdot\nabla T_k(w) + \int_\Omega j_s(x,w,\nabla w)T_k(w) > \sigma\|T_k(u)\|_{1,2}.$$

for every $w \in H_0^1(\Omega)$ with $\|w-u\|_{1,2} < \delta_1$.
Indeed, take $w_n \in H_0^1(\Omega)$ such that $w_n \to u$ in $H_0^1(\Omega)$. Then, (4.13) follows directly. Moreover, notice that from (1.3) and (1.4) it follows
$$j_s(x,w_n(x),\nabla w_n(x))w_n(x) \geqslant -R\beta(R)|\nabla w_n(x)|^2.$$

Since $w_n \to u$ in $H_0^1(\Omega)$, from (4.3) and by applying Fatou Lemma we get

$$\liminf_{n \to \infty} \left\{ \int_\Omega j_\xi(x, w_n, \nabla w_n) \cdot \nabla T_k(w_n) + \int_\Omega j_s(x, w_n, \nabla w_n) T_k(w_n) \right\}$$
$$\geq \int_\Omega j_\xi(x, u, \nabla u) \cdot \nabla T_k(u) + \int_\Omega j_s(x, u, \nabla u) T_k(u) > \sigma \|T_k(u)\|_{1,2}.$$

Thus, (4.14) holds. Let us now consider the continuous function $\mathcal{H} : B(u, \delta_1) \times [0, \delta_1] \to H_0^1(\Omega)$ defined by

$$\mathcal{H}(w, t) = w - \frac{t}{\|T_k(u)\|_{1,2}(1 + \varepsilon)} T_k(w).$$

From (4.13) and (4.14) we deduce that there exists $\delta < \delta_1$ such that

$$d(\mathcal{H}(w, t), w) \leq t$$
$$J(\mathcal{H}(w, t)) - J(w) \leq -\frac{\sigma}{1 + \varepsilon}$$

for every $t \in [0, \delta]$ and $w \in H_0^1(\Omega)$ with $\|w - u\|_{1,2} < \delta$ and $J(w) < J(u) + \delta$. Then, the arbitrariness of $\varepsilon$ yields $|dJ|(u) \geq \sigma$. Therefore, for every $k \geq 1$ we get

$$\int_\Omega j_s(x, u, \nabla u) T_k(u) + \int_\Omega j_\xi(x, u, \nabla u) \cdot \nabla T_k(u) \leq |dJ|(u) \|T_k(u)\|_{1,2}.$$

Taking the limit as $k \to +\infty$, by the Monotone Convergence Theorem one obtains (4.12). $\square$

Notice that a generalized solution $u$ (see definition 1.1) is not, in general, a distributional solution. This, because a test function $v \in W_u$ may not belong to $C_0^\infty$. Thus, it is natural to study the conditions under which it is possible to enlarge the class of admissible test functions. This kind of argument was introduced in [7]. More precisely, suppose we have a function $u \in H_0^1(\Omega)$ such that

(4.15) $$\int_\Omega j_\xi(x, u, \nabla u) \cdot \nabla z + \int_\Omega j_s(x, u, \nabla u) z = \langle w, z \rangle \qquad \forall z \in V_u,$$

where $V_u$ is defined in (4.4) and $w \in H^{-1}(\Omega)$. A natural question is whether or not we can take as test function $v \in H_0^1(\Omega) \cap L^\infty(\Omega)$. The next result gives an answer to this question.

**Theorem 4.7.** *Assume that conditions (1.1), (1.2), (1.3) hold. Let $w \in H^{-1}(\Omega)$ and $u \in H_0^1(\Omega)$ that satisfies (4.15). Moreover, suppose that $j_\xi(x, u, \nabla u) \cdot \nabla u \in L^1(\Omega)$ and there exist $v \in H_0^1(\Omega) \cap L^\infty(\Omega)$ and $\eta \in L^1(\Omega)$ such that*

(4.16) $$j_s(x, u, \nabla u) v + j_\xi(x, u, \nabla u) \cdot \nabla v \geq \eta.$$

*Then $j_\xi(x, u, \nabla u) \cdot \nabla v + j_s(x, u, \nabla u) v \in L^1(\Omega)$ and*

$$\int_\Omega j_\xi(x, u, \nabla u) \cdot \nabla v + \int_\Omega j_s(x, u, \nabla u) v = \langle w, v \rangle.$$

*Proof.* Since $v \in H_0^1(\Omega) \cap L^\infty(\Omega)$, then $H(\frac{u}{k}) v \in V_u$. From (4.15) we have

(4.17) $$\int_\Omega j_\xi(x, u, \nabla u) \cdot \nabla \left[ H\left(\frac{u}{k}\right) v \right] + \int_\Omega j_s(x, u, \nabla u) [H\left(\frac{u}{k}\right) v = \left\langle w, H\left(\frac{u}{k}\right) v \right\rangle$$

for every $k \geq 1$. Note that

$$\int_\Omega \left| j_\xi(x, u, \nabla u) \cdot \nabla u H'\left(\frac{u}{k}\right) \frac{v}{k} \right| \leq \frac{2}{k} \|v\|_\infty \int_\Omega j_\xi(x, u, \nabla u) \cdot \nabla u.$$

Since $j_\xi(x, u, \nabla u) \cdot \nabla u \in L^1(\Omega)$ we can use Lebesgue dominated convergence Theorem and deduce

$$\lim_{k \to \infty} \int_\Omega j_\xi(x, u, \nabla u) \cdot \nabla u H'\left(\frac{u}{k}\right)\frac{v}{k} = 0,$$

$$\lim_{k \to \infty} \left\langle w, H\left(\frac{u}{k}\right) v\right\rangle = \langle w, v\rangle.$$

As far as concerns the remaining terms in (4.17), notice that from (4.16) it follows

$$[j_s(x, u, \nabla u)v + j_\xi(x, u, \nabla u) \cdot \nabla v]H\left(\frac{u}{k}\right) \geq H\left(\frac{u}{k}\right)\eta \geq -\eta^- \in L^1(\Omega).$$

Thus, we can apply Fatou Lemma and obtain

$$\int_\Omega j_s(x, u, \nabla u)v + j_\xi(x, u, \nabla u) \cdot \nabla v \leq \langle w, v\rangle.$$

The previous inequality and (4.16) imply that

(4.18) $$j_s(x, u, \nabla u)v + j_\xi(x, u, \nabla u) \cdot \nabla v \in L^1(\Omega).$$

Now, notice that

$$\left|[j_s(x, u, \nabla u)v + j_\xi(x, u, \nabla u) \cdot \nabla v] H\left(\frac{u}{k}\right)\right| \leq |j_s(x, u, \nabla u)v + j_\xi(x, u, \nabla u) \cdot \nabla v|.$$

From (4.18) we deduce that we can use Lebesgue Dominated Convergence Theorem to pass to the limit in (4.17) and obtain the conclusion. □

In the next result we find the conditions under which we can use $v \in H_0^1(\Omega)$ in (4.15). Moreover, we prove, under suitable hypotheses, that if $u$ satisfies (4.15) then $u$ is a generalized solution (see Definition 1.1) of the corresponding problem.

**Theorem 4.8.** *Assume that conditions (1.1), (1.2), (1.3), (1.4) hold. Let $w \in H^{-1}(\Omega)$, and let $u \in H_0^1(\Omega)$ be such that (4.15) is satisfied. Moreover, suppose that $j_\xi(x, u, \nabla u) \cdot \nabla u \in L^1(\Omega)$, and that there exist $v \in H_0^1(\Omega)$ and $\eta \in L^1(\Omega)$ such that*

(4.19) $$j_s(x, u, \nabla u)v \geq \eta, \quad j_\xi(x, u, \nabla u) \cdot \nabla v \geq \eta.$$

*Then $j_s(x, u, \nabla u)v \in L^1(\Omega)$, $j_\xi(x, u, \nabla u) \cdot \nabla v \in L^1(\Omega)$ and*

(4.20) $$\int_\Omega j_\xi(x, u, \nabla u) \cdot \nabla v + \int_\Omega j_s(x, u, \nabla u)v = \langle w, v\rangle.$$

*In particular, it results $j_s(x, u, \nabla u)u, j_s(x, u, \nabla u) \in L^1(\Omega)$ and*

$$\int_\Omega j_\xi(x, u, \nabla u) \cdot \nabla u + \int_\Omega j_s(x, u, \nabla u)u = \langle w, u\rangle,$$

*moreover, $u$ is a generalized solution of the problem*

(4.21) $$\begin{cases} -\operatorname{div}(j_\xi(x, u, \nabla u)) + j_s(x, u, \nabla u) = w & \text{in } \Omega, \\ u = 0 & \text{on } \partial\Omega, \end{cases}$$

*in the sense of Definition 1.1.*

*Proof.* Let $k \geqslant 1$. For every $v \in H_0^1(\Omega)$, we have that $T_k(v) \in H_0^1(\Omega) \cap L^\infty(\Omega)$ and $-v^- \leqslant T_k(v) \leqslant v^+$, then from (4.19) we get

$$(4.22) \qquad j_s(x, u, \nabla u) T_k(v) \geqslant -\eta^- \in L^1(\Omega).$$

Moreover,

$$(4.23) \qquad j_\xi(x, u, \nabla u) \cdot \nabla T_k(v) \geqslant -[j_\xi(x, u, \nabla u) \cdot \nabla T_k(v)]^- \geqslant -\eta^- \in L^1(\Omega).$$

Then, we can apply Theorem 4.7 and obtain

$$(4.24) \qquad \int_\Omega j_s(x, u, \nabla u) T_k(v) + \int_\Omega j_\xi(x, u, \nabla u) \cdot \nabla T_k(v) = \langle w, T_k(v) \rangle$$

for every $k \geqslant 1$. By using again (4.22) and (4.23) and by arguing as in Theorem 4.7 we obtain

$$j_s(x, u, \nabla u) v \in L^1(\Omega), \quad j_\xi(x, u, \nabla u) \cdot \nabla v \in L^1(\Omega).$$

Thus, we can use Lebesgue Dominated Convergence Theorem to pass to the limit in (4.24) and obtain (4.20). In particular, by (1.4) and (4.3) we can choose $v = u$. Finally, since

$$j_s(x, u, \nabla u) = j_s(x, u, \nabla u) \chi_{\{|u|<1\}} + j_s(x, u, \nabla u) \chi_{\{|u| \geqslant 1\}}$$

and

$$\left| j_s(x, u, \nabla u) \chi_{\{|u| \geqslant 1\}} \right| \leqslant |j_s(x, u, \nabla u) u|,$$

by (1.3) it results also $j_s(x, u, \nabla u) \in L^1(\Omega)$. Finally, notice that if $v \in W_u$ we can take $\eta = j_\xi(x, u, \nabla u) \cdot \nabla v$ and $\eta = j_s(x, u, \nabla u) v$, so that (4.20) is satisfied. Thus, $u$ is a generalized solution of Problem (4.21). $\square$

Notice that the previous result implies that if $u \in H_0^1(\Omega)$ satisfies (4.15) and $j_\xi(x, u, \nabla u) \cdot \nabla u \in L^1(\Omega)$, it results that $j_s(x, u, \nabla u) \in L^1(\Omega)$, then $j_s(x, u, \nabla u) v \in L^1(\Omega)$ for every $v \in C_0^\infty(\Omega)$. While, the term that has not a distributional interpretation in (4.15) is $j_\xi(x, u, \nabla u)$. In the next result we show that if we multiply $j_\xi(x, u, \nabla u)$ by a suitable sequence of $C_c^1$ functions, we obtain, passing to the limit, a distributional interpretation of (4.15).

**Theorem 4.9.** *Assume conditions (1.1), (1.2), (1.3), (1.4). Let $w \in H^{-1}(\Omega)$ and $u \in H_0^1(\Omega)$ such that (4.15) is satisfied. Let $(\vartheta_h)$ be a sequence in $C_c^1(\mathbb{R})$ with*

$$\sup_h \|\vartheta_h\|_\infty < +\infty \qquad \sup_h \|\vartheta_h'\|_\infty < +\infty$$

$$\lim_{h \to +\infty} \vartheta_h(s) = 1, \qquad \lim_{h \to +\infty} \vartheta_h'(s) = 0.$$

*If $j_\xi(x, u, \nabla u) \cdot \nabla u \in L^1(\Omega)$, the sequence*

$$\mathrm{div}\left[\vartheta_h(u) j_\xi(x, u, \nabla u)\right]$$

*is strongly convergent in $W^{-1,q}(\Omega)$ for every $1 < q < \frac{N}{N-1}$ and*

$$\lim_{h \to +\infty} \left\{ -\mathrm{div}\left[\vartheta_h(u) j_\xi(x, u, \nabla u)\right] \right\} + j_s(x, u, \nabla u) = w \quad \text{in } W^{-1,q}(\Omega).$$

*Proof.* Let $w = -\operatorname{div} F$ with $F \in L^2(\Omega, \mathbb{R}^N)$ and $v \in C_c^\infty(\Omega)$. Then $\vartheta_h(u)v \in V_u$ and we can take $v$ as test function in (4.15). It results

$$\int_\Omega j_\xi(x, u, \nabla u)\vartheta_h(u)\nabla v = -\int_\Omega j_\xi(x, u, \nabla u)\vartheta_h'(u)\nabla u\, v - \int_\Omega j_s(x, u, \nabla u)\vartheta_h(u)v$$
$$+ \int_\Omega F\vartheta_h'(u)\nabla u\, v + \int_\Omega F\vartheta_h(u)\nabla v.$$

Then $u$ is a solution of the following equation

$$-\operatorname{div}\left[\vartheta_h(u)j_\xi(x, u, \nabla u)\right] = \xi_h \quad \text{in } \mathcal{D}'(\Omega),$$

where

$$\xi_h = -\left[\vartheta_h'(u)(j_\xi(x, u, \nabla u) - F)\cdot\nabla u + \vartheta_h(u)j_s(x, u, \nabla u)\right] - \operatorname{div}(\vartheta_h(u)F).$$

Now, notice that

$$\vartheta_h(u)F \to F, \qquad \text{strongly in } L^2(\Omega).$$

Then, $\operatorname{div}(\vartheta_h(u)F)$ is a convergent sequence in $H^{-1}(\Omega)$. Since the embedding of $H^{-1}(\Omega)$ in $W^{-1,q}(\Omega)$ is continuous, we get the desired convergence. Moreover, Theorem 4.8 implies that $j_s(x, u, \nabla u) \in L^1(\Omega)$. Then the remaining terms in $\xi_h$ converge strongly in $L^1(\Omega)$. Thus, we get the conclusion by observing that the embedding of $L^1(\Omega)$ in $W^{-1,q}(\Omega)$ is continuous. $\square$

Consider the case $j(x, s, \xi) = a(x, s)|\xi|^2$ with $a(x, s)$ measurable with respect to $x$, continuous with respect to $s$ and such that hypotheses (1.1), (1.2), (1.3), (1.4) are satisfied. The next result proves-in particular-that if there exists $u \in H_0^1(\Omega)$ that satisfies (4.15) and if $a(x, u)|\nabla u|^2 \in L^1(\Omega)$, then $u$ satisfies (4.15) in the sense of distribution.

**Theorem 4.10.** *Assume conditions (1.1), (1.2), (1.3), (1.4), (1.7). Let $w \in H^{-1}(\Omega)$ and $u \in H_0^1(\Omega)$ that satisfies (4.15). Moreover, suppose that $j_\xi(x, u, \nabla u)\cdot\nabla u \in L^1(\Omega)$ and that*

$$(4.25) \qquad\qquad\qquad j(x, s, \xi) = \widehat{j}(x, s, |\xi|).$$

*Then $j_\xi(x, u, \nabla u) \in L^1(\Omega)$ and $u$ is a solution of*

$$\begin{cases} -\operatorname{div}(j_\xi(x, u, \nabla u)) + j_s(x, u, \nabla u) = w & \text{in } \Omega, \\ u = 0 & \text{on } \partial\Omega, \end{cases}$$

*in $\mathcal{D}'(\Omega)$.*

*Proof.* It is readily seen that, in view of (1.1) and (4.25), it results

$$|\xi||j_\xi(x, s, \xi)| \leqslant \sqrt{2}j_\xi(x, s, \xi)\cdot\xi.$$

for a.e. $x \in \Omega$, every $s \in \mathbb{R}$ and $\xi \in \mathbb{R}^N$. Then

$$j_\xi(x, u, \nabla u)\chi_{\{|\nabla u|>1\}} \in L^1(\Omega).$$

Moreover, we take into account (1.7), and we observe that (1.5) implies that there exists a positive constant $C$ such that

$$|\xi| \leqslant 1 \implies |j_\xi(x, s, \xi)| \leqslant 4\alpha(|s|) \leqslant C(|s|^{p-2} + 1)$$

which by Sobolev embedding implies also that $j_\xi(x,u,\nabla u)\chi_{\{|\nabla u|\leqslant 1\}} \in L^1(\Omega)$. Then $j_\xi(x,u,\nabla u) \in L^1(\Omega)$. Moreover, from (1.3) and (1.4) we have

$$j_s(x,u,\nabla u)u \geqslant j_s(x,u,\nabla u)u\chi_{\{x:|u(x)|<R\}} \in L^1(\Omega).$$

Then Theorem 4.8 implies that $j_s(x,u,\nabla u)u \in L^1(\Omega)$. Finally, Theorem 4.8 yields the conclusion. □

## 5. A compactness Result for $J$

In this section we will prove the following compactness result for $J$. We will follow an argument similar to the one used in [9] and in [20].

**Theorem 5.1.** *Assume conditions (1.1), (1.2), (1.3), (1.4). Let $\{u_n\} \subset H_0^1(\Omega)$ be a bounded sequence with $j_\xi(x,u_n,\nabla u_n) \cdot \nabla u_n \in L^1(\Omega)$ and let $\{w_n\} \subset H^{-1}(\Omega)$ be such that*

$$(5.1) \qquad \forall v \in V_{u_n}: \int_\Omega j_s(x,u_n,\nabla u_n)v + j_\xi(x,u_n,\nabla u_n)\cdot \nabla v = \langle w_n, v\rangle.$$

*If $w_n$ is strongly convergent in $H^{-1}(\Omega)$, then, up to a subsequence, $u_n$ is strongly convergent in $H_0^1(\Omega)$.*

*Proof.* Let $w$ be the limit of $\{w_n\}$ and let $L > 0$ be such that

$$(5.2) \qquad \|u_n\|_{1,2} \leqslant L.$$

From (5.2) we deduce that there exists $u \in H_0^1(\Omega)$ such that, up to a subsequence,

$$(5.3) \qquad u_n \rightharpoonup u, \qquad \text{weakly in } H_0^1(\Omega).$$

*Step 1.* Let us first prove that $u$ is such that

$$(5.4) \qquad \int_\Omega j_\xi(x,u,\nabla u)\cdot \nabla \psi + \int_\Omega j_s(x,u,\nabla u)\psi = \langle w, \psi\rangle, \qquad \forall \psi \in V_u.$$

First of all, notice that from Rellich Compact Embedding Theorem we get that, up to a subsequence,

$$(5.5) \qquad \begin{cases} u_n \to u & \text{in } L^q(\Omega) \quad \forall q \in [1, 2N/(N-2)), \\ u_n(x) \to u(x) & \text{a.e. in } \Omega. \end{cases}$$

We now want to prove that, up to a subsequence,

$$(5.6) \qquad \nabla u_n(x) \to \nabla u(x) \quad \text{a.e. in } \Omega.$$

Let $h \geqslant 1$. For every $v \in C_c^\infty(\Omega)$ we have that $H\left(\frac{u_n}{h}\right)v \in V_{u_n}$ (where $H(s)$ is defined in (4.5)), then

$$\int_\Omega H\left(\frac{u_n}{h}\right) j_\xi(x,u_n,\nabla u_n) \cdot \nabla v$$
$$= -\int_\Omega \left[ H\left(\frac{u_n}{h}\right) j_s(x,u_n,\nabla u_n) + H'\left(\frac{u_n}{h}\right) j_\xi(x,u_n,\nabla u_n) \cdot \frac{\nabla u_n}{h}\right] v$$
$$+ \left\langle w_n, H\left(\frac{u_n}{h}\right) v\right\rangle.$$

Let $w_n = -\text{div}(F_n)$, with $(F_n)$ strongly convergent in $L^2(\Omega, \mathbb{R}^N)$. Then it follows that

$$\int_\Omega H\left(\frac{u_n}{h}\right) j_\xi(x, u_n, \nabla u_n) \cdot \nabla v$$
$$= \int_\Omega \left[ H'\left(\frac{u_n}{h}\right)(F_n - j_\xi(x, u_n, \nabla u_n)) \cdot \frac{\nabla u_n}{h} - H\left(\frac{u_n}{h}\right) j_s(x, u_n, \nabla u_n) \right] v$$
$$+ \int_\Omega H\left(\frac{u_n}{h}\right) F_n \cdot \nabla v.$$

Since the square bracket is bounded in $L^1(\Omega)$ and $(H\left(\frac{u_n}{h}\right) F_n)$ is strongly convergent in $L^2(\Omega, \mathbb{R}^N)$ we can apply [11, Theorem 5] with

$$b_n(x, \xi) = H\left(\frac{u_n(x)}{h}\right) j_\xi(x, u_n(x), \xi), \qquad E = E_h = \{x \in \Omega : |u(x)| \leqslant h\}$$

and deduce (5.6) by the arbitrariness of $h \geqslant 1$. Notice that by Theorem 4.8 we have for every $n$

$$\int_\Omega j_\xi(x, u_n, \nabla u_n) \cdot \nabla u_n + \int_\Omega j_s(x, u_n, \nabla u_n) u_n = \langle w_n, u_n \rangle.$$

Then, in view of (1.4) one has

(5.7) $$\sup_n \int_\Omega j_\xi(x, u_n, \nabla u_n) \cdot \nabla u_n < +\infty.$$

Let now $k \geqslant 1$, $\varphi \in C_c^\infty(\Omega)$, $\varphi \geqslant 0$ and consider

(5.8) $$v = \varphi e^{-M_k(u_n+R)^+} H\left(\frac{u_n}{k}\right), \quad M_k = \frac{\beta(2k)}{\alpha_0}.$$

Note that $v \in V_{u_n}$ and

$$\nabla v = \nabla \varphi e^{-M_k(u_n+R)^+} H\left(\frac{u_n}{k}\right) - M_k \varphi e^{-M_k(u_n+R)^+} \nabla(u_n+R)^+ H\left(\frac{u_n}{k}\right)$$
$$+ \varphi e^{-M_k(u_n+R)^+} H'\left(\frac{u_n}{k}\right) \frac{\nabla u_n}{k}.$$

When we take $v$ as test function in (5.1), we obtain

(5.9) $$\int_\Omega j_\xi(x, u_n, \nabla u_n) \cdot e^{-M_k(u_n+R)^+} H\left(\frac{u_n}{k}\right) \nabla \varphi$$
$$+ \int_\Omega \left[ j_s(x, u_n, \nabla u_n) - M_k j_\xi(x, u_n, \nabla u_n) \cdot \nabla(u_n+R)^+ \right] \varphi e^{-M_k(u_n+R)^+} H\left(\frac{u_n}{k}\right)$$
$$= \int_\Omega j_\xi(x, u_n, \nabla u_n) \cdot \varphi e^{-M_k(u_n+R)^+} H'\left(\frac{u_n}{k}\right) \frac{\nabla u_n}{k}$$
$$+ \left\langle w_n, \varphi e^{-M_k(u_n+R)^+} H\left(\frac{u_n}{k}\right) \right\rangle.$$

Observe that

$$\left[ j_s(x, u_n, \nabla u_n) - M_k j_\xi(x, u_n, \nabla u_n) \cdot \nabla(u_n+R)^+ \right] \varphi e^{-M_k(u_n+R)^+} H\left(\frac{u_n}{k}\right) \leqslant 0.$$

Indeed, the assertion follows from (1.4), for almost every $x$ such that $u_n(x) \leqslant -R$ while, for almost every $x$ in $\{x : -R \leqslant u_n(x) \leqslant 2k\}$ from (1.5), (4.3) and (5.8) we get

$$\left[ j_s(x, u_n, \nabla u_n) - M_k j_\xi(x, u_n, \nabla u_n) \cdot \nabla(u_n+R)^+ \right] \leqslant (\beta(2k) - \alpha_0 M_k)|\nabla u_n|^2 \leqslant 0.$$

Moreover, from (1.5), (5.2), (5.5) and (5.6) we have

$$\int_\Omega j_\xi(x, u_n, \nabla u_n) e^{-M_k(u_n+R)^+} H\left(\frac{u_n}{k}\right) \nabla \varphi \to \int_\Omega j_\xi(x, u, \nabla u) e^{-M_k(u+R)^+} H\left(\frac{u}{k}\right) \nabla \varphi$$

$$\left\langle w_n, \varphi e^{-M_k(u_n+R)^+} H\left(\frac{u_n}{k}\right) \right\rangle \to \left\langle w, \varphi e^{-M_k(u+R)^+} H\left(\frac{u}{k}\right) \right\rangle,$$

as $n \to \infty$. Moreover, we take into account (5.7) and deduce that there exists a positive constant $C > 0$ such that

$$\left| \int_\Omega j_\xi(x, u_n, \nabla u_n) \cdot \varphi e^{-M_k(u_n+R)^+} H'\left(\frac{u_n}{k}\right) \frac{\nabla u_n}{k} \right| \leqslant \frac{C}{k}.$$

We take the superior limit in (5.9) and we apply Fatou Lemma to obtain

$$(5.10) \quad \int_\Omega j_\xi(x, u, \nabla u) \cdot e^{-M_k(u+R)^+} H\left(\frac{u}{k}\right) \nabla \varphi + \int_\Omega j_s(x, u, \nabla u) \varphi e^{-M_k(u+R)^+} H\left(\frac{u}{k}\right)$$

$$- M_k \int_\Omega j_\xi(x, u, \nabla u) \cdot \nabla u^+ \varphi e^{-M_k(u+R)^+} H\left(\frac{u}{k}\right)$$

$$\geqslant -\frac{C}{k} + \left\langle w, \varphi e^{-M_k(u+R)^+} H\left(\frac{u}{k}\right) \right\rangle$$

for every $\varphi \in C_c^\infty(\Omega)$ with $\varphi \geqslant 0$. Then, the previous inequality holds for every $\varphi \in H_0^1 \cap L^\infty(\Omega)$ with $\varphi \geqslant 0$. We now choose in (5.10) the admissible test function

$$\varphi = e^{M_k(u+R)^+} \psi, \quad \psi \in V_u, \quad \psi \geqslant 0.$$

It results

$$(5.11) \quad \int_\Omega j_\xi(x, u, \nabla u) \cdot H\left(\frac{u}{k}\right) \nabla \psi + \int_\Omega j_s(x, u, \nabla u) H\left(\frac{u}{k}\right) \psi$$

$$\geqslant -\frac{C}{k} + \left\langle w, H\left(\frac{u}{k}\right) \psi \right\rangle.$$

Note that

$$\left| j_\xi(x, u, \nabla u) \cdot H\left(\frac{u}{k}\right) \nabla \psi \right| \leqslant |j_\xi(x, u, \nabla u)| |\nabla \psi|,$$

$$\left| j_s(x, u, \nabla u) H\left(\frac{u}{k}\right) \psi \right| \leqslant |j_s(x, u, \nabla u) \psi|.$$

Since $\psi \in V_u$ and from (1.3) and (1.5) we deduce that we can pass to the limit in (5.11) as $k \to +\infty$, and we obtain

$$\int_\Omega j_\xi(x, u, \nabla u) \cdot \nabla \psi + \int_\Omega j_s(x, u, \nabla u) \psi \geqslant \langle w, \psi \rangle, \quad \forall \psi \in V_u, \psi \geqslant 0.$$

In order to show the opposite inequality, we can take $v = \varphi e^{-M_k(u_n-R)^-} H\left(\frac{u_n}{k}\right)$ as test function in (5.1) and we can repeat the same argument as before. Thus, (5.4) follows.

<u>Step 2.</u>
In this step we will prove that $u_n \to u$ strongly in $H_0^1(\Omega)$. From (4.3), (5.7) and Fatou Lemma, we have

$$0 \leqslant \int_\Omega j_\xi(x, u, \nabla u) \cdot \nabla u \leqslant \liminf_n \int_\Omega j_\xi(x, u_n, \nabla u_n) \cdot \nabla u_n < +\infty$$

so that $j_\xi(x, u, \nabla u) \cdot \nabla u \in L^1(\Omega)$. Therefore, by Theorem 4.8 we deduce

$$(5.12) \qquad \int_\Omega j_\xi(x, u, \nabla u) \cdot \nabla u + \int_\Omega j_s(x, u, \nabla u) u = \langle w, u \rangle.$$

In order to prove that $u_n$ converges to $u$ strongly in $H_0^1(\Omega)$ we follow the argument of [20, Theorem 3.2] and we consider the function $\zeta : \mathbb{R} \to \mathbb{R}$ defined by

$$(5.13) \qquad \zeta(s) = \begin{cases} Ms & \text{if } 0 < s < R, \quad M = \frac{\beta(R)}{\alpha_0}, \\ MR & \text{if } s \geqslant R, \\ -Ms & \text{if } -R < s < 0, \\ MR & \text{if } s \leqslant -R. \end{cases}$$

We have that $v_n = u_n e^{\zeta(u_n)}$ belongs to $H_0^1(\Omega)$, and conditions (1.3), (1.4) and (1.5) imply that hypotheses of Theorem 4.8 are satisfied. Then, we deduce that we can use $v_n$ as test function in (5.1). It results

$$\int_\Omega j_\xi(x, u_n, \nabla u_n) \cdot \nabla u_n e^{\zeta(u_n)} = \langle w_n, v_n \rangle$$
$$- \int_\Omega \left[ j_s(x, u_n, \nabla u_n) + j_\xi(x, u_n, \nabla u_n) \cdot \nabla u_n \zeta'(u_n) \right] v_n$$

Note that $v_n$ converges to $u e^{\zeta(u)}$ weakly in $H_0^1(\Omega)$ and almost everywhere in $\Omega$. Moreover, conditions (1.3), (1.4) and (5.13) allow us to apply Fatou Lemma and get that

$$(5.14) \qquad \limsup_h \int_\Omega j_\xi(x, u_n, \nabla u_n) \cdot \nabla u_n e^{\zeta(u_n)} \leqslant \langle w, u e^{\zeta(u)} \rangle$$
$$- \int_\Omega \left[ j_s(x, u, \nabla u) + j_\xi(x, u, \nabla u) \cdot \nabla u \zeta'(u) \right] u e^{\zeta(u)}.$$

On the other hand (5.12) and (5.13) imply that

$$(5.15) \qquad \begin{cases} j_\xi(x, u, \nabla u) \cdot \nabla \left[ u e^{\zeta(u)} \right] + j_s(x, u, \nabla u) u e^{\zeta(u)} \in L^1(\Omega), \\ j_\xi(x, u, \nabla u) \cdot \nabla \left[ u\, e^{\zeta(u)} \right] \in L^1(\Omega). \end{cases}$$

Therefore, from Theorem 4.8 it results

$$(5.16) \qquad \int_\Omega j_\xi(x, u, \nabla u) \cdot \nabla (u e^{\zeta(u)}) + \int_\Omega j_s(x, u, \nabla u) u e^{\zeta(u)} = \langle w, u e^{\zeta(u)} \rangle.$$

Thus, (5.14) and (5.16) imply that

$$\int_\Omega j_\xi(x, u, \nabla u) \cdot \nabla u\, e^{\zeta(u)} \leqslant \liminf_{n \to \infty} \int_\Omega j_\xi(x, u_n, \nabla u_n) \cdot \nabla u_n e^{\zeta(u_n)}$$
$$\leqslant \limsup_{n \to \infty} \int_\Omega j_\xi(x, u_n, \nabla u_n) \cdot \nabla u_n e^{\zeta(u_n)} \leqslant \int_\Omega j_\xi(x, u, \nabla u) \cdot \nabla u e^{\zeta(u)}.$$

Then (4.3) implies that $u_n \to u$ strongly in $H_0^1(\Omega)$. □

## 6. Proofs of Theorems 2.1 and 2.3

In this section we give the definition of a Concrete Palais–Smale sequence, we study the relation between a Palais–Smale sequence and a Concrete Palais–Smale sequence, and we prove that $f$ satisfies the $(PS)_c$ for every $c \in \mathbb{R}$. Finally, we give the proofs of Theorems 2.1 and 2.3.

Let us consider the functional $I : H_0^1(\Omega) \to \mathbb{R}$ defined by

$$I(v) = -\int_\Omega G(x,v) - \langle \Lambda, v \rangle$$

where $\Lambda \in H^{-1}(\Omega)$, $G(x,s) = \int_0^s g(x,t)\,dt$ and $g : \Omega \times \mathbb{R} \to \mathbb{R}$ is a Carathéodory function satisfying assumption (2.2).

Then (1.2) implies that the functional $f : H_0^1(\Omega) \to \mathbb{R} \cup \{+\infty\}$ defined by $f(v) = J(v) + I(v)$ is lower semicontinuous by (1.2).

In order to apply the abstract theory the following result will be fundamental.

**Theorem 6.1.** *Assume conditions (1.1), (1.2), (1.4), (2.2). Then, for every $(u, \eta) \in \mathrm{epi}\,(f)$ with $f(u) < \eta$ it results*

$$|d\mathcal{G}_f|(u, \eta) = 1.$$

*Moreover, if $j(x, -s, -\xi) = j(x, s, \xi)$, $g(x, -s) = -g(x, s)$ and $\Lambda = 0$, for every $\eta > f(0)$ one has $|d_{\mathbb{Z}_2}\mathcal{G}_f|(0, \eta) = 1$.*

*Proof.* Since $G$ is of class $C^1$, Theorem 3.11 and Proposition 3.7 imply the result. □

Furthermore, since $G$ a $C^1$ functional, as a consequence of Proposition 4.5 one has the following result.

**Proposition 6.2.** *Assume conditions (1.2), (1.3), (1.5), (2.2). Let us consider $u \in \mathrm{dom}(f)$ with $|df|(u) < +\infty$. Then there exists $w \in H^{-1}(\Omega)$ such that $\|w\|_{-1,2} \leq |df|(u)$ and*

$$\forall v \in V_u : \int_\Omega j_\xi(x, u, \nabla u) \cdot \nabla v + \int_\Omega j_s(x, u, \nabla u)v - \int_\Omega g(x, u)v - \langle \Lambda, v \rangle = \langle w, v \rangle.$$

*Proof.* Given $u \in \mathrm{dom}(f)$ with $|df|(u) < +\infty$, let

$$\widehat{J}(v) = J(v) - \int_\Omega g(x, u)v - \langle \Lambda, v \rangle$$

$$\widehat{I}(v) = I(v) + \int_\Omega g(x, u)v.$$

Then, since $\widehat{I}$ is of class $C^1$ with $\widehat{I}'(u) = 0$, by (c) of Proposition 3.7 we get $|df|(u) = |d\widehat{J}|(u)$. By Proposition 4.5, there exists $w \in H^{-1}(\Omega)$ with $\|w\|_{-1,2} \leq |df|(u)$ and

$$\forall v \in V_u : \int_\Omega j_\xi(x, u, \nabla u) \cdot \nabla v + \int_\Omega j_s(x, u, \nabla u)v - \int_\Omega g(x, u)v - \langle \Lambda, v \rangle = \langle w, v \rangle$$

and the assertion is proved. □

We can now give the definition of the Concrete Palais–Smale condition.

**Definition 6.3.** *Let $c \in \mathbb{R}$. We say that $\{u_n\}$ is a Concrete Palais–Smale sequence for $f$ at level $c$ ($(CPS)_c$–sequence for short) if there exists $w_n \in H^{-1}(\Omega)$ with $w_n \to 0$ such that $j_\xi(x, u_n, \nabla u_n) \cdot \nabla u_n \in L^1(\Omega)$ for every $n$ and*

$$(6.1) \qquad f(u_n) \to c,$$

$$(6.2) \quad \int_\Omega j_\xi(x, u_n, \nabla u_n) \cdot \nabla v + \int_\Omega j_s(x, u_n, \nabla u_n) v - \int_\Omega g(x, u_n) v - \langle \Lambda, v \rangle$$
$$= \langle w_n, v \rangle, \quad \forall\, v \in V_{u_n}.$$

*We say that $f$ satisfies the Concrete Palais–Smale condition at level $c$ ($(CPS)_c$ for short) if every $(CPS)_c$–sequence for $f$ admits a subsequence strongly convergent in $H_0^1(\Omega)$.*

**Proposition 6.4.** *Assume conditions (1.2), (1.3), (1.5), (2.2). If $u \in \operatorname{dom}(f)$ satisfies $|df|(u) = 0$, then $u$ is a generalized solution of*

$$\begin{cases} -\operatorname{div}(j_\xi(x, u, \nabla u)) + j_s(x, u, \nabla u) = g(x, u) + \Lambda & \text{in } \Omega, \\ u = 0 & \text{on } \partial\Omega. \end{cases}$$

*Proof.* Combine Lemma 4.6, Proposition 6.2 and Theorem 4.8. □

The following result concerns the relation between the $(PS)_c$ condition and the $(CPS)_c$ condition.

**Proposition 6.5.** *Assume conditions (1.2), (1.3), (1.5), (2.2). Then if $f$ satisfies the $(CPS)_c$ condition, it satisfies the $(PS)_c$ condition.*

*Proof.* Let $\{u_n\} \subset \operatorname{dom}(f)$ that satisfies the Definition 3.4. From Lemma 4.6 and Proposition 6.2 we get that $u_n$ satisfies the conditions in Definition 6.3. Thus, there exists a subsequence, which converges in $H_0^1(\Omega)$. □

We now want to prove that $f$ satisfies the $(CPS)_c$ condition at every level $c$. In order to do this, let us consider a $(CPS)_c$–sequence $\{u_n\} \in \operatorname{dom}(f)$.

From Theorem 5.1 we deduce the following result.

**Proposition 6.6.** *Assume conditions (1.1), (1.2), (1.3), (1.4), (2.2). Let $\{u_n\}$ be a $(CPS)_c$–sequence for $f$, bounded in $H_0^1(\Omega)$. Then $\{u_n\}$ admits a strongly convergent subsequence in $H_0^1(\Omega)$.*

*Proof.* Let $\{u_n\} \subset \operatorname{dom}(f)$ be a concrete Palais–Smale sequence for $f$ at level $c$. Taking into account that, as known, by (2.2) the map $\{u \mapsto g(x, u)\}$ is compact from $H_0^1(\Omega)$ to $H^{-1}(\Omega)$, it suffices to apply Theorem 5.1 to see that $\{u_n\}$ is strongly compact in $H_0^1(\Omega)$. □

**Proposition 6.7.** *Assume conditions (1.1), (1.2), (1.3), (1.4), (1.9), (2.2), (2.3). Then every $(CPS)_c$–sequence $\{u_n\}$ for $f$ is bounded in $H_0^1(\Omega)$.*

*Proof.* Condition (1.4) and (4.3) allow us to apply Theorem 4.8 to deduce that we may choose $v = u_n$ as test functions in (6.2). Taking into account conditions (1.9), (2.2), (2.7), (6.1), the boundedness of $\{u_n\}$ in $H_0^1(\Omega)$ follows by arguing as in [20, Lemma 4.3]. □

**Remark 6.8.** *Notice that we use condition (1.9) only in Proposition 6.7.*

Under the assumptions of Theorem 2.3 we now have the following result.

**Theorem 6.9.** *Assume conditions* (1.1), (1.2), (1.3), (1.4), (1.9), (2.8), (2.3). *Then the functional $f$ satisfies the $(PS)_c$ condition at every level $c \in \mathbb{R}$.*

*Proof.* Let $\{u_n\} \subset \text{dom}(f)$ be a concrete Palais–Smale sequence for $f$ at level $c$. From Proposition 6.7 it follows that $\{u_n\}$ is bounded in $H_0^1(\Omega)$. By Proposition 6.6 $f$ satisfies the Concrete Palais–Smale condition. Finally Proposition 6.5 implies that $f$ satisfies the $(PS)_c$ condition. □

We are now able to prove Theorem 2.1.

*Proof of Theorem 2.1.*
We will prove Theorem 2.1 as a consequence of Theorem 3.10. First, note that (1.2) and (2.2) imply that $f$ is lower semicontinuous. Moreover, from (2.5) we deduce that $f$ is an even functional, and from Theorem 3.11 we deduce that (3.3) and condition (d) of Theorem 3.10 are satisfied. Hypotheses (1.7) and (2.4) imply that condition (b) of Theorem 3.10 is verified.

Let now $(\lambda_h, \varphi_h)$ be the sequence of solutions of $-\Delta u = \lambda u$ with homogeneous Dirichlet boundary conditions. Moreover, let us consider $V^+ = \overline{\text{span}}\{\varphi_h \in H_0^1(\Omega) : h \geqslant h_0\}$ and note that $V^+$ has finite codimension. In order to prove (a) of Theorem 3.10 it is enough to show that there exist $h_0, \gamma > 0$ such that

$$\forall u \in V^+ : \|\nabla u\|_2 = 1 \implies f(u) \geqslant \gamma.$$

First, note that condition (2.2) implies that for every $\varepsilon > 0$ we find $a_\varepsilon^{(1)} \in C_c^\infty(\Omega)$ and $a_\varepsilon^{(2)} \in L^{\frac{2N}{N+2}}(\Omega)$ with $\|a_\varepsilon^{(2)}\|_{\frac{2N}{N+2}} \leqslant \varepsilon$ and

$$|g(x,s)| \leqslant a_\varepsilon^{(1)}(x) + a_\varepsilon^{(2)}(x) + \varepsilon|s|^{\frac{N+2}{N-2}}.$$

Now, let $u \in V^+$ and notice that there exist $c_1, c_2 > 0$ such that

$$\begin{aligned}
f(u) &\geqslant \alpha_0 \|\nabla u\|_2^2 - \int_\Omega G(x,u) \\
&\geqslant \alpha_0 \|\nabla u\|_2^2 - \int_\Omega \left(\left(a_\varepsilon^{(1)} + a_\varepsilon^{(2)}\right)|u| + \frac{N-2}{2N}\varepsilon|u|^{\frac{2N}{N-2}}\right) \\
&\geqslant \alpha_0 \|\nabla u\|_2^2 - \|a_\varepsilon^{(1)}\|_2 \|u\|_2 - c_1 \|a_\varepsilon^{(2)}\|_{\frac{2N}{N+2}} \|\nabla u\|_2 - \varepsilon c_2 \|\nabla u\|_2^{\frac{2N}{N-2}} \\
&\geqslant \alpha_0 \|\nabla u\|_2^2 - \|a_\varepsilon^{(1)}\|_2 \|u\|_2 - c_1 \varepsilon \|\nabla u\|_2 - \varepsilon c_2 \|\nabla u\|_2^{\frac{2N}{N-2}}.
\end{aligned}$$

Then if $h_0$ is sufficiently large, since $\lambda_h \to +\infty$, for all $u \in V^+$, $\|\nabla u\|_2 = 1$ implies $\|a_\varepsilon^{(1)}\|_2 \|u\|_2 \leqslant \alpha_0/2$. Thus, for $\varepsilon > 0$ small enough, $\|\nabla u\|_2 = 1$ implies $f(u) \geqslant \gamma$ for some $\gamma > 0$. Then also (a) of Theorem 3.10 is satisfied.

Theorem 6.9 implies that $f$ satisfies $(PS)_c$ condition at every level $c$, so that we get the existence of a sequence of critical points $\{u_h\} \subset H_0^1(\Omega)$ with $f(u_h) \to +\infty$. Finally, Proposition 6.4 yields the conclusion. □

Let us conclude this section by proving Theorem 2.3.

*Proof of Theorem 2.3.*
We will prove Theorem 2.3 as a consequence of Theorem 3.9. In order to do this, let us notice that from (1.2) and (2.8) $f$ is lower semicontinuous on $H_0^1(\Omega)$. Moreover, Theorem 3.11 implies that condition (3.3) is satisfied. From Theorem 6.9 we deduce that $f$ satisfies $(PS)_c$ condition at every level $c$. It is left to show that $f$ satisfies the geometrical assumptions of Theorem (3.9).

Let us first consider the case in which $\Lambda = 0$, then we will apply the first part of Theorem 3.9. Notice that condition (1.2), (2.8) and (2.9) imply that there exist $\eta$ and $r$ that satisfy 3.4. Conditions (1.2) and (2.4) imply

$$(6.3) \qquad f(v) \leqslant \int_\Omega \alpha(|v|)|\nabla v|^2 - \int_\Omega k(x)|v|^p + \|a_0\|_1 + C_0 \|a_1\|_{\frac{2N}{N+2}} \|v\|_{1,2}.$$

Now, let us consider $W$ a finite dimensional subspace of $H_0^1(\Omega)$ such that $W \subset L^\infty(\Omega)$. Condition (1.7) implies that for every $\varepsilon > 0$ there exists $R > r$, $w \in W$, with $\|w\|_\infty > R$ and a positive constant $C_\varepsilon$ such that

$$(6.4) \qquad \int_\Omega \alpha(|w|)|\nabla w|^2 \leqslant \varepsilon C_W \|w\|_{1,2}^p + C_\varepsilon \|w\|_{1,2}^2,$$

where $C_W$ is a positive constant depending on $W$. Then (6.3) and (6.4) yield (3.5), by choosing suitably $\varepsilon$. Thus, we can apply Theorem 3.9 and deduce the existence of a nontrivial critical point $u \in H_0^1(\Omega)$. From Proposition 6.4 we get that $u$ is a generalized solution of Problem $(P)$.

Now, let us consider the case in which $\Lambda \not\equiv 0$. As before, we note that (1.2), (2.8) and (2.9) imply that there exists $\varepsilon > 0$ such that for every $\Lambda \in H^{-1}(\Omega)$ with $\|\Lambda\| \leqslant \varepsilon$, it results

$$f(u) \geqslant \gamma, \qquad \text{for every } u \text{ with } \|u\|_{1,2} = r.$$

Moreover, we use condition (1.2), (1.7) and (2.4) and we argue as before to deduce the existence of $v_1 \in H_0^1(\Omega)$ with $\|v_1\| > r$ and $f(v_1) < 0$. Finally, let us consider $\varphi_1$ the first eigenfunction of the Laplacian operator with Dirichlet homogeneous boundary conditions. We set $v_0 = t_0 \varphi_1$, with $t_0$ sufficiently small, so that, thanks to (1.2) and (2.8), also $f(v_0) < 0$. Then we can apply Theorem 3.9 and we get the existence of two critical points $u_1, u_2 \in H_0^1(\Omega)$. Finally, Proposition 6.4 yields the conclusion. □

**Remark 6.10.** *Notice that Theorems 1.2 and 1.3 are an easy consequence of Theorems 2.1 and 2.3 respectively. Indeed, consider for example $g_1(x,s) = a(x)\mathrm{arctg}s + |s|^{p-2}s$, in order to prove Theorem 1.2, it is left to show that $g_1(x,s)$ satisfies conditions (2.2), (2.3) and (2.4). First, notice that Young inequality implies that, for every $\varepsilon > 0$ there exists a positive constant $\beta(\varepsilon)$ such that (2.2) holds with $a_\varepsilon(x) = \beta(\varepsilon) + a(x)$. Moreover, (2.3) is satisfied with $a_0(x) = 0$ and $b_0(x) = \pi/2(p-1)$. Finally, (2.4) is verified with $k(x) = 1/p$, $\overline{a}(x) = 0$ and $\overline{b}(x) = (\pi/2 + C)a(x)$ where $C \in \mathbb{R}^+$ is sufficiently large. Theorem 1.3 can be obtained as a consequence of Theorem 2.3 in an analogous way.*

## 7. Summability Results

In this section we suppose that $g(x,s)$ satisfies the following growth condition

$$(7.1) \qquad |g(x,s)| \leqslant a(x) + b|s|^{\frac{N+2}{N-2}}, \qquad a(x) \in L^r(\Omega),\ b \in \mathbb{R}^+.$$

Note that (2.2) implies (7.1). In this section we will set $2^* = 2N/(N-2)$. We will prove the following theorem.

**Theorem 7.1.** *Let us assume that conditions (1.1), (1.2), (1.3), (1.4), (2.2) are satisfied. Let $u \in H_0^1(\Omega)$ be a generalized solution of problem $(P)$. Then the following conclusions hold.*
*(a) If $r \in (2N/(N+2), N/2)$, $u$ belongs to $L^{r^{**}}(\Omega)$, where $r^{**} = Nr/(N-2r)$.*
*(b) If $r > N/2$, (where $r$ is defined in (2.8)) $u$ belongs to $L^\infty(\Omega)$.*

Theorem 7.1 will be proved as a consequence of the following lemma.

**Lemma 7.2.** *Let us assume that conditions (1.1), (1.2), (1.3), (1.4) are satisfied. Let $u \in H_0^1(\Omega)$ be a generalized solution of the problem*

(7.2) $$\begin{cases} -\operatorname{div}(j_\xi(x,u,\nabla u)) + j_s(x,u,\nabla u) + c(x)u = f(x) & \text{in } \Omega, \\ u = 0 & \text{on } \partial\Omega, \end{cases}$$

*then the following conclusions hold.*

(a') *If $c \in L^{\frac{N}{2}}(\Omega)$ and $f \in L^r(\Omega)$, with $r \in (2N/(N+2), N/2)$, $u$ belongs to $L^{r^{**}}(\Omega)$, where $r^{**} = Nr/(N-2r)$.*

(b') *If $c \in L^t(\Omega)$ with $t > N/2$ and $f \in L^q(\Omega)$, with $q > N/2$, $u$ belongs to $L^\infty(\Omega)$.*

*Proof.* Let us first prove conclusion a'). For every $k > R$ (where $R$ is defined in (1.4)), let us define the function $\eta_k(s) : \mathbb{R} \to \mathbb{R}$ such that $\eta_k \in C^1$, $\eta_k$ is odd and

(7.3) $$\eta_k(s) = \begin{cases} 0 & \text{if } 0 < s < R, \\ (s-R)^{2\gamma+1} & \text{if } R < s < k, \\ b_k s + c_k & \text{if } s > k, \end{cases}$$

where $b_k$ and $c_k$ are constant such that $\eta_k$ is $C^1$. Since $u$ is a generalized solution of (7.2), $v = \eta_k(u)$ belongs to $W_u$. Then we can take it as test function, moreover, $j_s(x,u,\nabla u)\eta_k(u) \geqslant 0$. Then from (1.4) and (4.3) we get

(7.4) $$\alpha_0 \int_\Omega \eta_k'(u)|\nabla u|^2 \leqslant \int_\Omega f(x)\eta_k(u) + \int_\Omega c(x)u\eta_k(u).$$

Now, let us consider the odd function $\psi_k(s) : \mathbb{R} \to \mathbb{R}$ defined by

(7.5) $$\psi_k(s) = \int_0^s \sqrt{\eta_k'(t)}\, dt.$$

The following properties of the functions $\psi_k$ and $\eta_k$ can be deduced from (7.3) and (7.5) by easy calculations

(7.6) $$\left[\psi_k'(s)\right]^2 = \eta_k'(s),$$

(7.7) $$0 \leqslant \eta_k(s)s \leqslant C_0 \psi_k(s)^2,$$

(7.8) $$|\eta_k(s)| \leqslant C_0 |\psi_k(s)|^{\frac{2\gamma+1}{\gamma+1}},$$

where $C_0$ is a positive constant. Notice that for every $\varepsilon > 0$ there exist $c_1(x) \in L^{\frac{N}{2}}(\Omega)$, with $\|c_1\|_{\frac{N}{2}} \leqslant \varepsilon$ and $c_2 \in L^\infty(\Omega)$ such that $c(x) = c_1(x) + c_2(x)$. From (7.4), (7.6), (7.7) and Hölder inequality, we deduce

$$\alpha_0 \int_\Omega |\nabla(\psi_k(u))|^2 \leqslant C_0 \|c_1(x)\|_{\frac{N}{2}} \left[\int_\Omega |\psi_k(u)|^{2^*}\right]^{\frac{2}{2^*}} + \int_\Omega |f(x) + C_0 c_2(x)u||\eta_k(u)|.$$

We fix $\varepsilon = (\alpha_0 \mathcal{S})/(2C_0)$, where $\mathcal{S}$ is the Sobolev constant. We obtain

(7.9) $$\left[\int_\Omega |\psi_k(u)|^{2^*}\right]^{\frac{2}{2^*}} \leqslant C \int_\Omega |f(x) + C_0 c_2(x)u||\eta_k(u)|.$$

Now, let us define the function

(7.10) $$h(x) = |f(x) + C_0 c_2(x)u(x)|,$$

and note that $h(x)$ belongs to $L^t(\Omega)$ with

(7.11) $$t = \min\{r, 2^*\}.$$

Let us consider first the case in which $t = r$, then from (7.8) and (7.9), we get

$$\left[\int_\Omega |\psi_k(u)|^{2^*}\right]^{\frac{2}{2^*}} \leqslant C\|h\|_r \left[\int_\Omega |\psi_k(u)|^{r'\frac{2\gamma+1}{\gamma+1}}\right]^{\frac{1}{r'}}.$$

Since $2N/(N+2) < r < N/2$ we can define $\gamma \in \mathbb{R}^+$ by

(7.12) $$\gamma = \frac{r(N+2) - 2N}{2(N-2r)} \implies 2^*(\gamma+1) = r'(2\gamma+1) = r^{**}$$

Moreover, since $r < N/2$ we have that $2/2^* > 1/r'$, then

(7.13) $$\left[\int_\Omega |\psi_k(u)|^{2^*}\right]^{\frac{2}{2^*} - \frac{1}{r'}} \leqslant C\|h\|_r.$$

Note that $|\psi_k(u)|^{2^*} \to C(\gamma)|u - R|^{\gamma+1}\chi_{\{|u(x)|>R\}}$ almost everywhere in $\Omega$. Then Fatou Lemma implies that $|u - R|^{\gamma+1}\chi_{\{|u(x)|>R\}}$ belongs to $L^{2^*}(\Omega)$. Thus, $u$ belongs to $L^{2^*(\gamma+1)}(\Omega) = L^{r^{**}}(\Omega)$ and the conclusion follows. Consider now the case in which $t = 2^*$ and note that this implies that $N > 6$. In this case we get

$$\left[\int_\Omega |\psi_k(u)|^{2^*}\right]^{\frac{2}{2^*}} \leqslant C\|h\|_{2^*} \left[\int_\Omega s|\psi_k(u)|^{(2^*)'\frac{2\gamma+1}{\gamma+1}}\right]^{\frac{1}{(2^*)'}}.$$

Since $N > 6$ it results $2/2^* > 1/(2^*)'$. Moreover, we can choose $\gamma$ such that

$$2^*(\gamma+1) = (2^*)'(2\gamma+1).$$

Thus, we follow the same argument as in the previous case and we deduce that $u$ belongs to $L^{s_1}(\Omega)$ where

$$s_1 = \frac{2^*N}{N - 2\,2^*}.$$

If it still holds $s_1 < r$ we can repeat the same argument to gain more summability on $u$. In this way for every $s \in [2^*, r)$ we can define the increasing sequence

$$s_0 = 2^*, \qquad s_{n+1} = \frac{Ns_n}{N - 2s_n},$$

and we deduce that there exists $\overline{n}$ such that $s_{\overline{n}-1} < r$ and $s_{\overline{n}} \geqslant r$. At this step from (7.11) we get that $t = r$ and then $u \in L^{r^{**}}(\Omega)$, that is the maximal summability we can achieve.

Now, let us prove conclusion $b'$). First, note that since $f \in L^q(\Omega)$, with $q > N/2$, $f$ belongs to $L^r(\Omega)$ for every $r > (2N)/(N+2)$. Then, conclusion $a'$) implies that $u \in L^\sigma(\Omega)$ for every $\sigma > 1$. Now, take $\delta > 0$ such that $t - \delta > N/2$, since $u \in L^{\frac{t}{\delta}}(\Omega)$ it results

$$\int_\Omega |c(x)u(x)|^{t-\delta} \leqslant \|c(x)\|_t^{t-\delta} \left[\int_\Omega |u(x)|^{\frac{t}{\delta}}\right]^{\frac{\delta}{t}} < \infty.$$

Then, the function $d(x) = f(x) - c(x)u(x)$ belongs to $L^r(\Omega)$ with $r = \min\{q, t - \delta\} > N/2$. Let us take $k > R$ ($R$ is defined in (1.4)) and consider the function $v = G_k(u) = u - T_k(u)$ (where $T_k(s)$ is defined in (3.6)). Since $u$ is a generalized solution of (7.2) we can take $v$ as test function. From (1.4) and (4.3) it results

$$\alpha_0 \int_\Omega |\nabla G_k(u)|^2 \leqslant \int_\Omega |d(x)||G_k(u)|.$$

The conclusion follows from Theorem 4.2 of [21]. □

**Remark 7.3.** *In the classical results of this type (see [18] or [6]) it is usually used as test function $v = |u|^{2\gamma}u$. Note that this type of function cannot be used here as test function because it does not belong to $W_u$. Moreover, the classical truncation $T_u()$ seem not to be useful because of the presence of $c(x)u$. Then we have chosen a suitable truncation of $u$ in order to manage the term $c(x)u$.*

Now we are able to prove Theorem 7.1.

*Proof of Theorem 7.1.*

Theorem 7.1 will be proved as a consequence of Lemma 7.2. So, consider $u$ a generalized solution of Problem $(P)$, we have to prove that $u$ is a generalized solution of Problem 7.2 for suitable $f(x)$ and $c(x)$. This is shown in Theorem 2.2.5 of [9], then we will give here a sketch of the proof of [9] just for clearness.
We set

$$g_0(x, s) = \min\{\max\{g(x, s), -a(x)\}, a(x)\}, \quad g_1(x, s) = g(x, s) - g_0(x, s).$$

It follows that $g(x, s) = g_0(x, s) + g_1(x, s)$ and $|g_0(x, s)| \leqslant a(x)$ so that we can set $f(x) = g_0(x, u(x))$. Moreover, we define

$$c(x) = \begin{cases} -\dfrac{g_1(x, u(x))}{u(x)} & \text{if } u(x) \neq 0, \\ 0 & \text{if } u(x) = 0. \end{cases}$$

Then $|c(x)| \leqslant b|u(x)|^{\frac{4}{N-2}}$, so that $c(x) \in L^{\frac{N}{2}}(\Omega)$. Lemma 7.2 implies that conclusion a) holds. Now, if $r > N/2$ we have that $f(x) \in L^r(\Omega)$ with $r > N/2$. Moreover, conclusion a) implies that $u \in L^t(\Omega)$ for every $t < \infty$, so that $c(x) \in L^t(\Omega)$ with $t > N/2$. Then Lemma 7.2 implies that $u \in L^\infty(\Omega)$. □

**Remark 7.4.** *When dealing with quasilinear equations (i.e. $j(x, s, \xi) = a(x, s)\xi \cdot \xi$), a standard technique, to prove summability results, is to reduce the problem to the linear one and to apply the classical result due to [21]. Note that here this is not possible because of the general form of $j_\xi(x, s, \xi)$.*

**Acknowledgments.** The authors wish to thank Marco Degiovanni for providing helpful discussions.

Dipartimento di Matematica
P.le Aldo Moro 2, I–00185 Roma, Italy.
*E-mail address*: pellacci@mat.uniroma1.it

Dipartimento di Matematica e Fisica
Via Musei 41, I–25121 Brescia, Italy.
*E-mail address*: m.squassina@dmf.unicatt.it